\documentclass[11pt,letter,reqno]{amsart}
\usepackage[latin9]{inputenc}
\usepackage{amstext}
\usepackage{amsthm}
\usepackage{amssymb}
\usepackage{graphicx}
\usepackage{comment}
\usepackage{amsmath}
\usepackage{breqn}
\makeatletter
\numberwithin{equation}{section}
\numberwithin{figure}{section} 
\theoremstyle{plain}

 \theoremstyle{remark}

\newtheorem{theorem}{Theorem}[section]

\newtheorem{lemma}[theorem]{Lemma}

\newtheorem{remark}[theorem]{Remark}

\usepackage{amsthm}
\usepackage{color}
\usepackage{color}\usepackage{xcolor}
\usepackage{subcaption}
\usepackage{enumerate}
\usepackage[draft]{todonotes}
\usepackage[american]{babel}
\usepackage{geometry}

\usepackage{color,soul}

\newcommand{\ri}{\mathrm{i}}

\newcommand{\diag}{\mathrm{diag}}
\newcommand{\tr}{\mathrm{trace}}
\newcommand{\erf}{\mathrm{erf}}
\newcommand{\sign}{\mathrm{sign}}
\newcommand{\erfi}{\mathrm{erfi}}
\newcommand{\adj}{\mathrm{adj}}

\def\Id {{ \rm Id}}

\makeatother

\providecommand{\remarkname}{Remark}
\providecommand{\theoremname}{Theorem}

\begin{document}

\title[]{Spectral Closure for the Linear Boltzmann-BGK Equation}

\author{Florian Kogelbauer}
\address{ETH Z\"{u}rich, Department of Mechanical and Process Engineering, Leonhardstrasse 27, 8092 Z\"{u}rich, Switzerland}
\email{floriank@ethz.ch}

\author{Ilya Karlin}
\address{ETH Z\"{u}rich, Department of Mechanical and Process Engineering, Leonhardstrasse 27, 8092 Z\"{u}rich, Switzerland}
\email{ikarlin@ethz.ch}

\begin{abstract}
We give an explicit description of the spectral closure for the three-dimensional linear Boltzmann-BGK equation in terms of the macroscopic fields, density, flow velocity and temperature. This results in a new linear fluid dynamics model which is valid for any relaxation time. The non-local exact fluid dynamics equations are compared to the Euler, Navier--Stokes and Burnett equations. Our results are based on a detailed spectral analysis of the linearized Boltzmann-BGK operator together with a suitable choice of spectral projection.   
\end{abstract}

\maketitle

\section{Introduction}


Since the invention of kinetic theory by Boltzmann \cite{Boltzmann1872} and Maxwell \cite{Maxwell1860}, the fundamental question arose: What is the connection between kinetic equations and the equations for the motion of continua?  Or, to phrase it differently: Can the governing equations of fluid dynamics be rigorously derived from kinetic theory?\\
%

This problem has a long history. Famously, in his speech at the International  Congress of Mathematics in Paris in 1900, Hilbert proposed a program to derive the passage from the atomistic view of fluids and gases to the motion of continua \cite{hilbert2000mathematical}. One interpretation of this challenge, known as Hilbert's sixth problem in this context, aims to prove the convergence of kinetic models, such as the Boltzmann equation, to known hydrodynamic models such as the Euler and the Navier--Stokes equations \cite{SRay2007,saint2009hydrodynamic,saint2014mathematical}.\\

The derivation of hydrodynamics from kinetic models is often regarded as a closure problem where one seeks a self-consistent expression of the fluxes in the balance equations for primitive (macroscopic) fields of mass density, momentum density and energy density. On a formal level, a well-established approach to the closure problem is the Chapman--Enskog expansion \cite{Enskog1917,Chapman1916,chapman1990mathematical}, where a Taylor series in powers of the Knudsen number (the molecular mean free path to a characteristic flow scale ratio) is invoked. The lower-order approximations lead to compressible Euler and Navier--Stokes--Fourier systems. The undeniable success of the Chapman--Enskog method is rooted in the evaluation of the phenomenological transport coefficients, viscosity and thermal conductivity for a one-component gas, in terms of the microscopic interaction potential between particles, as well as prediction of the thermodiffusion effect in gas mixtures \cite{Brush1976}. On the other hand, extension of the Chapman--Enskog approximation beyond the classical Navier--Stokes--Fourier order, the Burnett and super-Burnett approximations \cite{Burnett1936,chapman1990mathematical}, encountered difficulties. Even in the simplest regime, while linearized around a global equilibrium, the higher-order hydrodynamic closure may exhibit an instability, as first shown by Bobylev \cite{bobylev1982chapman} for the Burnett and the super-Burnett approximations for Maxwell's molecules. Since the global equilibrium is stable by way of the dissipative nature of the Boltzmann equation, Bobylev's instability is an artifact brought about by the Chapman--Enskog procedure. The problem of higher-order hydrodynamics is exacerbated  in the non-linear regime. Indeed, as pointed out by Slemrod \cite{slemrod2012chapman}, convergence of a singular expansion to the leading-order equation is by no means obvious: the formation of shocks might be an obstacle to global uniform convergence in the sense of solutions \cite{SLEMROD20131497}. Furthermore, the expansion of a non-local operator in frequency space in terms of (local) differential operators may be problematic. As a remedy, Rosenau suggested a non-local closure \cite{PhysRevA.40.7193} based on rational functions rather than polynomial approximations to the Chapman--Enskog solution.\\


A different approach is to address the problem of hydrodynamics from kinetics as a problem of \emph{invariant manifolds}.  This viewpoint was first suggested in a short paper by McKean \cite{McKean1969} and expanded in a series of works by Gorban \& Karlin \cite{gorban1994method,gorban2005invariant,gorban2014hilbert}.
For model systems (Grad's moment systems \cite{grad1949kinetic}), it was shown that the method of invariant manifold is equivalent to \emph{exact summation} of the Chapman--Enskog series to all orders \cite{gorban1996short}.
We term the latter case \emph{exact hydrodynamics} since, once achieved, it furnishes the complete characterization of the hydrodynamic limit of the kinetic equation and hence, the rigorous and exact closure. In this setting, the problem remains non-trivial even in the linear case, for infinite-dimensional problems. Accurate numerical solutions were found in \cite{karlin2008exact}, on the level of the linear Boltzmann-BGK kinetic model \cite{bhatnagar1954model}, and extended to a finite-moment approximation of the linear Boltzmann equation for Maxwell's molecules in \cite{colangeli2009}.\\

In this work, the derivation of the exact (valid to all admissible scales) linear hydrodynamics is considered in two consecutive steps.
First, the slow invariant manifold is identified as the linear subspace spanned over the hydrodynamic spectrum of the linearized Boltzmann-BGK operator. This is achieved on the basis of the explicit solution to the eigenvalue problem presented recently in Kogelbauer \& Karlin \cite{kogelbauer2023exact}. Let us refer to \cite{grad1963asymptotic,ellis1975first,palczewski1984} for qualitative results on the spectra of general linear kinetic operators, including the existence of hydrodynamic branches, critical wave number local expansions for small wave numbers. In \cite{kogelbauer2020slow, kogelbauer2021, kogelbauer2023spectral}, explicit spectral calculations have been performed for several kinetic models, including explict expressions of critical wave numbers and branch merging.\\
However, the knowledge of the spectrum alone is not the final step towards the derivation of hydrodynamic equations. The next step is the projection of the dynamics onto the slow manifold in terms of primitive variables, density, momentum and energy (or temperature). To that end, we derive the hydrodynamic projection in two independent ways. First, we demonstrate that all information about the projection is essentially encoded in a function of eigenvalues, which we call spectral temperature. This direct computation uses specific features of the BGK model. On the other hand, the hydrodynamic projection can be equivalently derived on the basis of the Riesz spectral projector, a more general route applicable to a variety of linear kinetic problems. Both approaches are shown to be consistent with one another, and resulting in the unique hydrodynamic projection. Let us emphasize that we derive a closed-form expression for the transport coefficients in wave space (transport operators in physical space) in terms of eigenvalues.\\

The structure of the paper is as follows: Preliminaries in Sec.\ \ref{sec:notation} include the notation and nomenclature, in particular, the plasma dispersion function is introduced. Some useful properties of the plasma dispersion function necessary for the spectral analysis of the linearized Boltzmann-BGK operator are collected in Appendix \ref{propZ} for the sake of completeness. In Sec.\ \ref{sec:closure}, following the invariant manifold formulation of the closure problem \cite{gorban2005invariant}, we introduce the closure operator for a generic linear kinetic equation. While the majority of derivations of hydrodynamics proceed in terms of primitive variables to solve the invariance equation, our approach is different. We first recognize the slow invariant manifold from the analysis of the spectral problem and, secondly, induce the dynamics on this manifold in terms of primitive variables by a coordinate change from spectral variables to hydrodynamic fields. The realization of this program starts in Sec.\ \ref{sec:spectrum} where we first review analytical results on the spectral problem of the linearized Boltzmann-BGK model \cite{kogelbauer2023exact}.  Sec.\ \ref{sec:spectral_closure} is devoted to the explicit construction of the coordinate change from spectral variables to macroscopic variables, involving a single analytic function depending on eigenvalues, called \textit{spectral temperature}. In Sec. \ref{HydroEqu}, we present the exact hydrodynamic equations for each wave vector in frequency space, as well as in physical space. While, classically, the closure is obtained through transport coefficients relating the dynamics of the macroscopic variables to each other, the exact hydrodynamic equations involve transport operators with finite frequency support acting on the corresponding variables. Finally, in Sec.\ \ref{compare}, we compare the exact non-local hydrodynamics to local approximations such as the Euler equation, the Navier--Stokes--Fourier system and the Burnett system. In particular, we recover the approximate slow dynamics obtained through the Chapman--Enskog expansion. We conclude with a discussion in Sec.\ \ref{sec:conclusion}.


\section{Notation and Basic Definitions}
\label{sec:notation}
For a wave vector $\mathbf{k}\in\mathbb{Z}^3$, $\mathbf{k}=(k_1,k_2,k_3)$, we denote its wave number as 
\begin{equation}
k:=|\mathbf{k}|=\sqrt{k_1^2+k_2^2+k_3^2}.
\end{equation}
For a given wave vector $\mathbf{k}\neq 0$, we define a coordinate system with a component parallel and with two components orthogonal to $\mathbf{k}$ by splitting any vector $\mathbf{v}\in\mathbb{R}^3$ as
\begin{equation}\label{split}
    \mathbf{v}= \mathbf{v}_{\parallel}+\mathbf{v}_{\perp},
\end{equation}
where $\mathbf{v}_{\parallel}=\frac{1}{k^2}(\mathbf{v}\cdot\mathbf{k})\mathbf{k}$ and $ \mathbf{v}_{\perp}=-\frac{1}{k^2}\mathbf{k}\times(\mathbf{k}\times\mathbf{v})$, which satisfies $\mathbf{v}^{\perp}\cdot\mathbf{k} =0$. This can be achieved by a rotation matrix $\mathbf{Q}_{\mathbf{k}}$ satisfying $\mathbf{k}=\mathbf{Q}_{\mathbf{k}}(k,0,0)^T$ to give 
\begin{equation}
    \mathbf{v}=\mathbf{Q}_{\mathbf{k}}(v_{\parallel},v_{\perp 1},v_{\perp 2}),
\end{equation}
where $v_{\parallel} = \mathbf{v}\cdot\mathbf{k} $ and $(v_{\perp 1},v_{\perp 2})$ are the components of the unit base vectors of $\mathbf{v}_{\perp}$. The matrix $\mathbf{Q}_{\mathbf{k}}$ can be determined by, e.g. the Rodrigues' rotation formula \cite{rodrigues1840lois}:
\begin{equation}\label{defQ}
    \mathbf{Q}_{\mathbf{k}} =\left(
\begin{array}{ccc}
 \frac{k_1}{k} & -\frac{k_2}{k} & -\frac{k_3}{k} \\
 \frac{k_2}{k} & 1-\frac{k_2^2}{k^2+k_1 k} & -\frac{k_2 k_3}{k^2+k_1 k} \\
 \frac{k_3}{k} & -\frac{k_2 k_3}{k^2+k_1 k} & 1-\frac{k_3^2}{k^2+k_1 k} \\
\end{array}
\right). 
\end{equation}
For later calculations, we also define the $5\times 5$ block-diagonal matrix
\begin{equation}\label{defQtilde}
    \tilde{\mathbf{Q}}_{\mathbf{k}} = \diag(1,\mathbf{Q}_{\mathbf{k}},1). 
\end{equation}

We introduce the \textit{plasma dispersion function} as the integral 
\begin{equation}\label{defZ}
    Z(\zeta) = \frac{1}{\sqrt{2\pi}} \int_{\mathbb{R}} \frac{e^{-\frac{v^2}{2}}}{v-\zeta}\, dv,
\end{equation}
for any $\zeta\in\mathbb{C}\setminus\mathbb{R}$. The function $Z$ is analytic on each half plane $\{\Im(\zeta)>0\}$ and $\{\Im(\zeta)<0\}$ and satisfies the complex differential equation
\begin{equation}\label{diffZ}
    \frac{dZ}{d\zeta} = -\zeta Z-1.
\end{equation}
As the name suggest, function \eqref{defZ} appears in plasma physics in the context of Landau damping \cite{fitzpatrick2014plasma}. We collect further useful properties of the plasma dispersion function in the Appendix \ref{propZ}.\\

Let $\mathcal{H}$ denote a Hilbert space and let $\mathbf{T}:\mathcal{H}\to\mathcal{H}$ be a linear operator with domain of definition $\mathcal{D}(\mathcal{H})$. We denote the spectrum of $\mathbf{T}$ as $\sigma(\mathbf{T})$ and its resolvent set as $\rho(\mathbf{T})$.
The main operator $\mathcal{L}_{\mathbf{k}}$ of this paper (to be defined later) will be defined on the Hilbert space
\begin{equation}
\mathcal{H}_{\mathbf{v}} = L^2_{\mathbf{v}}(\mathbb{R}^3,e^{-|\mathbf{v}|^2}),
\end{equation}
together with the inner product
\begin{equation}\label{innerprod}
\langle f, g \rangle_{\mathbf{v}} = (2\pi)^{-\frac{3}{2}}\int_{\mathbb{R}^3} f(\mathbf{v}) g^*(\mathbf{v}) e^{-\frac{|\mathbf{v}|^2}{2}} d\mathbf{v}.
\end{equation}
For later calculation and to ease notation, we define the following set of basis vectors:
\begin{equation}\label{base5}
\begin{split}
e_0(\mathbf{v}) &= (2\pi)^{-\frac{3}{4}},\\
e_1(\mathbf{v}) &= (2\pi)^{-\frac{3}{4}} v_1,\\
e_2(\mathbf{v}) &= (2\pi)^{-\frac{3}{4}} v_2,\\
e_3(\mathbf{v}) &= (2\pi)^{-\frac{3}{4}} v_3,\\
e_4(\mathbf{v}) &= (2\pi)^{-\frac{3}{4}} \frac{|\mathbf{v}|^2-3}{\sqrt{6}},
\end{split}
\end{equation} 
which satisfy the orthonormality relation
\begin{equation}\label{orthorel}
\langle e_i, e_j \rangle_{\mathbf{v}} = \delta_{ij},\quad \text{ for } \quad  0 \leq i,j \leq 4,
\end{equation}
where $\delta_{ij}$ is the Kronecker's delta. We bundle the basis functions \eqref{base5} into a single vector
\begin{equation}
    \mathbf{e}=(e_0,e_1,e_2,e_3,e_4). 
\end{equation}
For a one-body distribution function $f: \mathbb{T}^3\times\mathbb{R}^3\times [0,\infty)\to\mathbb{R}^{+}$, we introduce the moments of the distribution function $f$ as
\begin{equation}\label{defmomentum}
\mathbf{M}^{(n)}(\mathbf{x},t)=\int_{\mathbb{R}^3}f(\mathbf{x},\mathbf{v},t)\,\mathbf{v}^{\otimes n}d\mathbf{v},
\end{equation}
where $\mathbf{v}^{\otimes 0}=1$, $\mathbf{v}^{\otimes 1}=\mathbf{v}$ and 
\begin{equation}
\mathbf{v}^{\otimes n}=\underbrace{\mathbf{v}\otimes...\otimes \mathbf{v}}_{n-\text{times}},
\end{equation}
for $n\geq 2$ is the $n$-th tensor power. The moment defined in \eqref{defmomentum} is an $n$-th order symmetric tensor, depending on space and time.\\
Given a vector $X=(x_1,\dots,x_n)$, we denote the set of cyclical permutations of $X$ as
\begin{equation}\label{defcirc}
    \circlearrowright (x_1,\dots,x_n) = \{(x_1,x_2,\dots,x_n),(x_2,x_3,\dots,x_n,x_1),\dots,(x_n,x_1,\dots,x_{n-1})\}. 
\end{equation}
For a matrix $A$, we denote its adjugate as $\adj(A)$, which satisfies $A \,\adj(A)=\det(A)\Id$. \\
We denote the strip between $-a$ and $0$ as
\begin{equation}\label{defR}
    \mathcal{R}_a = \{z\in\mathbb{C}: -a< \Re z< 0\}\subset\mathbb{C}.
\end{equation}

\section{The Closure Problem for the Linear BGK Equation}
\label{sec:closure}

In this section, we recall the classical closure problem for kinetic equations in general and illustrate it on the BGK equation in particular. First, we formulate the governing equations suitable for our setting and illustrate the closure problems for the hierarchy of moment equations. Subsequently, we define the closure operator and outline the relation of the existence of a (slow) invariant manifold with an exact closure relation.\\

We will be interested in the three-dimensional BGK equation linearized around a global Maxwellian:
\begin{equation}\label{maineq}
\frac{\partial f}{\partial t}+\mathbf{v}\cdot\nabla_{\mathbf{x}} f=-\frac{1}{\tau}L_{\rm BGK}[f],
\end{equation}
for the deviation relative to the global Maxwellian $f: \mathbb{T}^3\times\mathbb{R}^3\times [0,\infty)\to\mathbb{R}$, $f=f(\mathbf{x},\mathbf{v},t)$ and the BGK collision operator
\begin{equation}
L_{\rm BGK}[f](\mathbf{x},\mathbf{v},t)=\Big(f(\mathbf{x},\mathbf{v},t)-\mathbb{P}_5[f](\mathbf{x},\mathbf{v},t)\Big).
\end{equation}
The projection operator $\mathbb{P}_5:\mathcal{H}_{\mathbf{v}}\to\mathcal{H}_{\mathbf{v}}$ is defined as
\begin{equation}\label{defP5}
\mathbb{P}_5f = \sum_{j=0}^4 \langle f, e_j \rangle_{\mathbf{v}} e_j,
\end{equation}
i.e., the projection onto the first five basis vectors \eqref{base5}. Clearly, \eqref{defP5} defines an orthogonal projection with respect to \eqref{innerprod}. Integrating equation \eqref{maineq} in $\mathbf{x}$ shows that the five basis functions \eqref{base5} are center modes (and the dynamic in these directions is conserved), since
\begin{equation}\label{center}
L_{\rm BGK}[e_j]=0, \quad \text{ for } \quad  0 \leq j \leq 4.
\end{equation}

Expanding $f$ in a Fourier series
\begin{equation}
f(\mathbf{x},\mathbf{v})= \sum_{|\mathbf{k}|=0}^{\infty}\hat{f}(\mathbf{k},\mathbf{v}) e^{\ri \mathbf{x}\cdot\mathbf{k}}, 
\end{equation}
for the Fourier coefficients
\begin{equation}
\hat{f}(\mathbf{k},\mathbf{v}) = \frac{1}{(2\pi)^3}\int_{\mathbb{R}^3} f(\mathbf{x},\mathbf{v})e^{-\ri \mathbf{x}\cdot\mathbf{k}}\, d\mathbf{x},
\end{equation}
the linear operator in \eqref{maineq} can be unitarily conjugated to the family of operators
\begin{equation}\label{defLhat}
    \mathcal{L}_{\mathbf{k}} =  -\ri\mathbf{v}\cdot\mathbf{k} - \frac{1}{\tau}(1-\mathbb{P}_5),
\end{equation}
indexed by the wave vector $\mathbf{k}$. 

Because of the normalization of the basis functions \eqref{base5}, the relation to the \textit{macroscopic variables} density $\rho$, velocity $\mathbf{u}$ and temperature $T$ is given by
\begin{equation}\label{macro}
\begin{split}
    \rho & = \langle f, e_0\rangle_{\mathbf{v}} = (2\pi)^{-\frac{3}{2}}\int_{\mathbb{R}^3} f(\mathbf{v}) e^{-\frac{|\mathbf{v}|^2}{2}} d\mathbf{v},\\
    \mathbf{u} & = \langle f, (e_1,e_2,e_3)\rangle_{\mathbf{v}} = (2\pi)^{-\frac{3}{2}}\int_{\mathbb{R}^3} f(\mathbf{v})\mathbf{v} e^{-\frac{|\mathbf{v}|^2}{2}} d\mathbf{v},\\
    T & =\sqrt{\frac{2}{3}} \langle f, e_4\rangle_{\mathbf{v}} = (2\pi)^{-\frac{3}{2}}\int_{\mathbb{R}^3} f(\mathbf{v})\frac{|\mathbf{v}|^2-3}{3} e^{-\frac{|\mathbf{v}|^2}{2}} d\mathbf{v},
\end{split}
\end{equation}
which we bundle into a single vector
\begin{equation}\label{defh} 
    \mathbf{h} = \left(\begin{array}{c}
    \rho\\
    \mathbf{u}\\
    \sqrt{\frac{3}{2}}T\end{array}\right).
\end{equation}
Because of the orthonormality relations \eqref{orthorel}, we prefer to work with the basis $(e_0,...,e_4)$. To account for the prefactor $\sqrt{{2}/{3}}$ in the definition of the temperature in \eqref{macro} in the final (physically meaningful) dynamical equations, we have to multiply the last entries accordingly.

\subsection{The Closure Problem}
Let us recall the classical closure problem for kinetic equations illustrated on the BGK equation.\\
Multiplying the BGK equation \eqref{maineq} with $\mathbf{v}^{\otimes n}$ and integrating in velocities gives the following hierarchy of moment equations
\begin{equation}\label{momentumeq}
    \frac{\partial}{\partial t} \mathbf{M}^{(n)} = - \nabla\cdot\mathbf{M}^{(n+1)} -\frac{1}{\tau}\mathbf{M}^{(n)} + \frac{1}{\tau}\mathbf{M}^{(n)}_{\rm eq,lin},
\end{equation}
where
\begin{equation}
    \mathbf{M}^{(n)}_{\rm eq,lin} = \mathbf{M}^{(0)}\langle \mathbf{v}^{\otimes n}, 1 \rangle + \langle \mathbf{v}^{\otimes n}, \mathbf{v}\cdot \mathbf{M}^{(1)}  \rangle +  \frac{1}{3}\tr(\mathbf{M}^{(2)}-\text{Id}_{3\times 3})\langle \mathbf{v}^{\otimes n}, \frac{|\mathbf{v}|^2-3}{3}\rangle.
\end{equation}
Ideally, we would like to obtain a closed system for the first few moments, thus allowing for a consistent dynamical description of macroscopic variables. As equation \eqref{momentumeq} illustrates, however, will the rate of change of $\mathbf{M}^{(n)}$ always be affected by the flux of the next moment $\nabla\cdot\mathbf{M}^{(n+1)}$ (the term $\mathbf{M}^{(n)}_{eq,lin}$ only comprises moments up to order two). Consequently, there is no way to obtain a self-consistent moment system from the full dynamics of the kinetic model \eqref{maineq}.\\
As a way out of this inconvenient matter of facts, we can constrain the dynamics of our system \eqref{maineq} by assuming that the full dynamics is given parametrically as a function of, say, the five macroscopic variables
\begin{equation}\label{closure1}
    f(\mathbf{x},\mathbf{v},t) = F(\rho(\mathbf{x},t),\mathbf{u}(\mathbf{x},t),T(\mathbf{x},t);\mathbf{v}).
\end{equation}
Writing equation \eqref{maineq} a bit more abstractly as 
\begin{equation}\label{maineqL}
    \frac{\partial f}{\partial t} = \mathcal{L}[f],\quad \mathcal{L} = -\mathbf{v}\cdot\nabla_{\mathbf{x}}-\frac{1}{\tau} L_{\rm BGK},
\end{equation}
and denoting $\mathbb{P}_5^{\perp} = 1-\mathbb{P}_5$, assumption \eqref{closure1} corresponds to the existence of a linear operator $\mathcal{C}:\text{range }\mathbb{P}_5\to \text{range }\mathbb{P}_{5}^{\perp}$, called \textit{closure operator}, such that 
\begin{equation}\label{oneplusC}
    f=(1+\mathcal{C})\mathbb{P}_5 f,
\end{equation}
and the dynamics of the macroscopic variables $\mathbb{P}_5 f$ can be written self-consistently as
\begin{equation}\label{dynh}
 \frac{\partial \mathbb{P}_5 f}{\partial t} = \mathbb{P}_5\mathcal{L}(1+\mathcal{C})\mathbb{P}_5 f,
\end{equation}
while the closure operator $\mathcal{C}$ satisfies the condition of \textit{exact closure}:
\begin{equation}\label{exactclosure}
(\mathcal{C}\mathbb{P}_5-\mathbb{P}_5^{\perp})\mathcal{L}(1+\mathcal{C})= 0. 
\end{equation}
Indeed, applying $\mathbb{P}_{5}^{\perp}$ to equation \eqref{maineqL} and using assumption \eqref{oneplusC}, we obtain
\begin{equation}
\begin{split}
 \mathbb{P}_5^{\perp}\mathcal{L}(1+\mathcal{C})\mathbb{P}_5 f = \frac{\partial}{\partial t} \mathbb{P}_5^{\perp} f = \mathcal{C}\frac{\partial\mathbb{P}_5 f}{\partial t}.
  \end{split}
\end{equation}
Using now the reduced dynamics \eqref{dynh}, we arrive at
\begin{equation}
    \mathbb{P}_5^{\perp}\mathcal{L}(1+\mathcal{C})\mathbb{P}_5 f = \mathcal{C}\mathbb{P}_5\mathcal{L}(1+\mathcal{C})\mathbb{P}_5 f,
\end{equation}
which is equivalent to \eqref{exactclosure}.
Since the operator $\mathcal{L}$ can be written as the direct sum over operators $\mathcal{L}_{\mathbf{k}}$ for $\mathbf{k}\in\mathbb{Z}^3$, we effectively seek a closure operator for each wave vector by writing $\hat{f}_{\mathbf{k}} = (1+\mathcal{C}_{\mathbf{k}})\mathbb{P}_5 \hat{f}_{\mathbf{k}}$ with 
\begin{equation}\label{dynhk}
 \frac{\partial \mathbb{P}_5 \hat{f}_{\mathbf{k}}}{\partial t} = \mathbb{P}_5\mathcal{L}_{\mathbf{k}}(1+\mathcal{C}_{\mathbf{k}})\mathbb{P}_5 \hat{f}_{\mathbf{k}},
\end{equation}
and the condition of being an exact closure at each wave vector:
\begin{equation}\label{exactclosurek}
(\mathcal{C}_{\mathbf{k}}\mathbb{P}_5-\mathbb{P}_5^{\perp})\mathcal{L}_{\mathbf{k}}(1+\mathcal{C}_{\mathbf{k}})= 0,
\end{equation}
for $\mathbf{k}\in\mathbb{Z}^3$.
\begin{remark}
Equation \eqref{exactclosurek} for the closure operator is a special case of the \emph{invariance equation}, and is the cornerstone of any derivation of hydrodynamics from kinetic theory. For example, the Chapman--Enskog method is based on a Taylor series expansion in terms of a small parameter $\epsilon$ after a rescaling $\tau \to \epsilon\tau$ while the method of invariant manifold \cite{gorban1994method} uses Newton-type iterations. Note that, even in the simplest linear setting addressed here, the invariance equation \eqref{exactclosurek} is non-linear (quadratic) in the unknown closure operator. In \cite{karlin2008exact}, a numerical solution to \eqref{exactclosurek} was obtained for the linear Boltzmann-BGK kinetic model.	
\end{remark}

\begin{remark}
Below, we shall circumvent solving the invariance equation \eqref{exactclosurek} directly. Instead, the exact closure operator shall be determined in two steps. First, we identify the slow invariant manifold based on the properties of the spectrum of $\mathcal{L}_{\mathbf{k}}$ (for details on the spectrum of the BGK equation, we refer to the following section). Indeed, for a certain range of wave numbers $0<k<k_{\rm crit,min}$, we will see that the spectral properties of $\mathcal{L}_{\mathbf{k}}$ allow us to define a closure operator $\mathcal{C}_{\mathbf{k}}$ which has the property that a general solution (restricted to $\mathbf{k}$) approaches the dynamics of the $\mathcal{C}_{\mathbf{k}}$-constrained ensemble exponentially fast in time. Second, we shall find a unique projection of the dynamics onto the slow invariant manifold in terms of primitive variables. The latter step is the main aspect of this work.
\end{remark}

\section{Spectrum of the Linear BGK Equation and Spectral Closure}
\label{sec:spectrum}

In this section, we first recall the properties of the BGK spectrum derived in \cite{kogelbauer2023exact}. This involves the three families of modal branches (diffusion, acoustics and shear), as well as their asymptotic behavior for small wave number. Then we define the hydrodynamic manifold as the eigenspace associated to the hydrodynamic modes and define the spectral closure. This will be achieved by a change of coordinates from spectral variables to macroscopic variables, thus providing a solution to \eqref{exactclosurek}.\\
Because the hydrodynamic modes are slow, any trajectory of distribution functions (deviations relative to the global Maxwellian in the linear case) will approach this linear manifold exponentially fast in time. Consequently, the general moment dynamics \eqref{momentumeq} will be approximate exponentially well in time with the self-consistent moment system derived from the spectral closure operator. 

\subsection{Properties of the Spectrum}

In order to define the spectral closure for the BGK system, we recall the most important implications of the detailed spectral analysis performed in \cite{kogelbauer2023exact}, including a complete description of the eigenvalues above the essential spectrum for each wave number as zeros of a holomorphic spectral function, as well as the Taylor expansion in wave number. We define the hydrodynamic manifold as a linear combination of eigenvectors and derive the spectral dynamics on the manifold.\\
\begin{theorem}
The spectrum of the linearized BGK operator $\mathcal{L}$ with relaxation time $\tau$ around a global Maxwellian is given by
\begin{equation}\label{thmspec}
\sigma(\mathcal{L}) = \left\{-\frac{1}{\tau}+\ri\mathbb{R}\right\}\cup\bigcup_{N\in \text{Modes}}\bigcup_{|\mathbf{k}|<k_{{\rm crit},N}}\{\lambda_{N}(\tau|\mathbf{k}|)\},
\end{equation}
where $\text{Modes}=\{\text{shear}, \text{diff}, \text{ac}, \text{ac}*\}$ corresponding to the shear mode, the diffusion mode and the pair of complex conjugate acoustic modes. The essential spectrum is given by the line $\Re\lambda=-\frac{1}{\tau}$, while the discrete spectrum consists of a \textit{finite} number of discrete, isolated eigenvalues. Along with each family of modes, there exists a critical wave number $k_{crit,N}$, limiting the range of wave numbers for which $\lambda_N$ exists. The modes $\{\text{diff}, \text{ac}, \text{ac}*\}$ all have algebraic multiplicity one, while the shear mode has algebraic and geometric multiplicity two. 
The eigenvalues (in dependence of the wave number $k$ and the relaxation time $\tau$) are given as zeros of the spectral function:
\begin{equation}\label{defSigma}
\begin{split}
&    \Sigma_{k,\tau}(\lambda)  = \frac{1}{6(\ri k\tau)^5}(Z(\zeta)-\ri\tau k)^2\\
&\times\Big(\zeta+6 \ri k^3 \tau ^3-\zeta  (\zeta^2+5) k^2 \tau ^2+2 \ri (\zeta ^2+3)k \tau -4 \ri Z^2 (\zeta )((\zeta ^2+1) k \tau -\ri \zeta)\\
&\qquad +Z(\zeta ) (\zeta ^2-(\zeta^4+4 \zeta ^2+11) k^2 \tau ^2+2 \ri  k \tau\zeta ^3 -5) ) \Big)\Big|_{\zeta = \ri\frac{\tau\lambda+1}{k\tau}}.
\end{split}
\end{equation}
\end{theorem}

\begin{figure}
	\begin{subfigure}{.45\textwidth}
		\centering
		\includegraphics[width=1\linewidth]{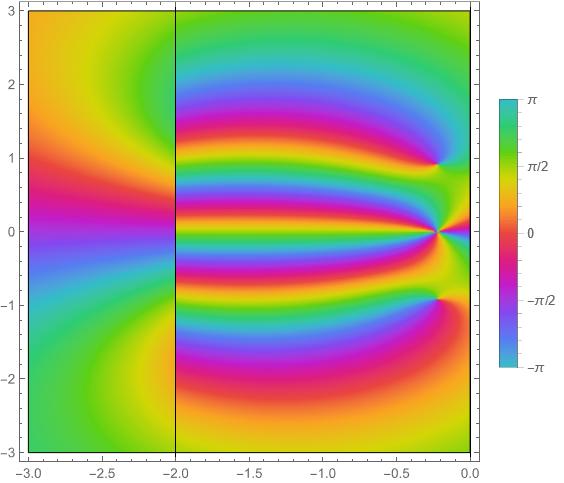}
		\caption{$k=0.7$}
		\label{p1}
	\end{subfigure}
	\begin{subfigure}{.45\textwidth}
		\centering
		\includegraphics[width=1\linewidth]{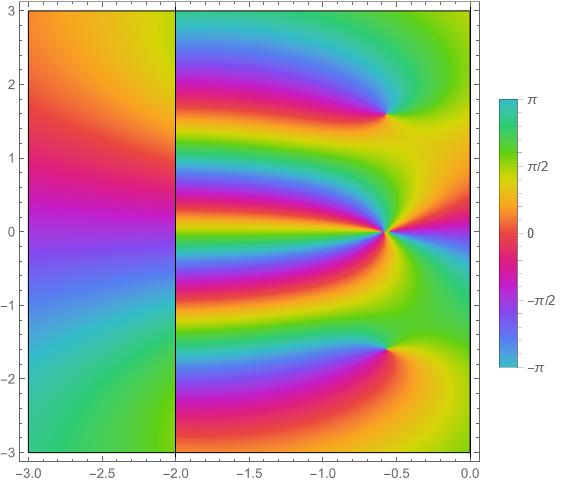}
		\caption{$k=1.2$}
		\label{ps2}
	\end{subfigure}
	\caption{Argument plot of the spectral function \eqref{defSigma} for $\tau=0.5$ and different values of $k$. The zeros of the function $\Sigma_{k,\tau}$ in the complex plane define eigenvalues of the linearized BGK operator.}
	\label{argplotSigma}
\end{figure}

For a proof, we refer to \cite{kogelbauer2023exact}.
A typical argument plot of spectral function \eqref{defSigma} is shown in Figure \ref{argplotSigma}. For $k=0$, the function \eqref{defSigma} collapses to a multiple of $\lambda^5$, recovering the center spectrum (conserved quantities) of \eqref{maineq}, see \eqref{center}. Increasing $k$, the zeros of  $\Sigma_{k,\tau}$ branch out and decrease monotonically in their real parts. 
\begin{figure}
    \begin{subfigure}{.5\textwidth}
		\centering
		\includegraphics[width=0.8\linewidth]{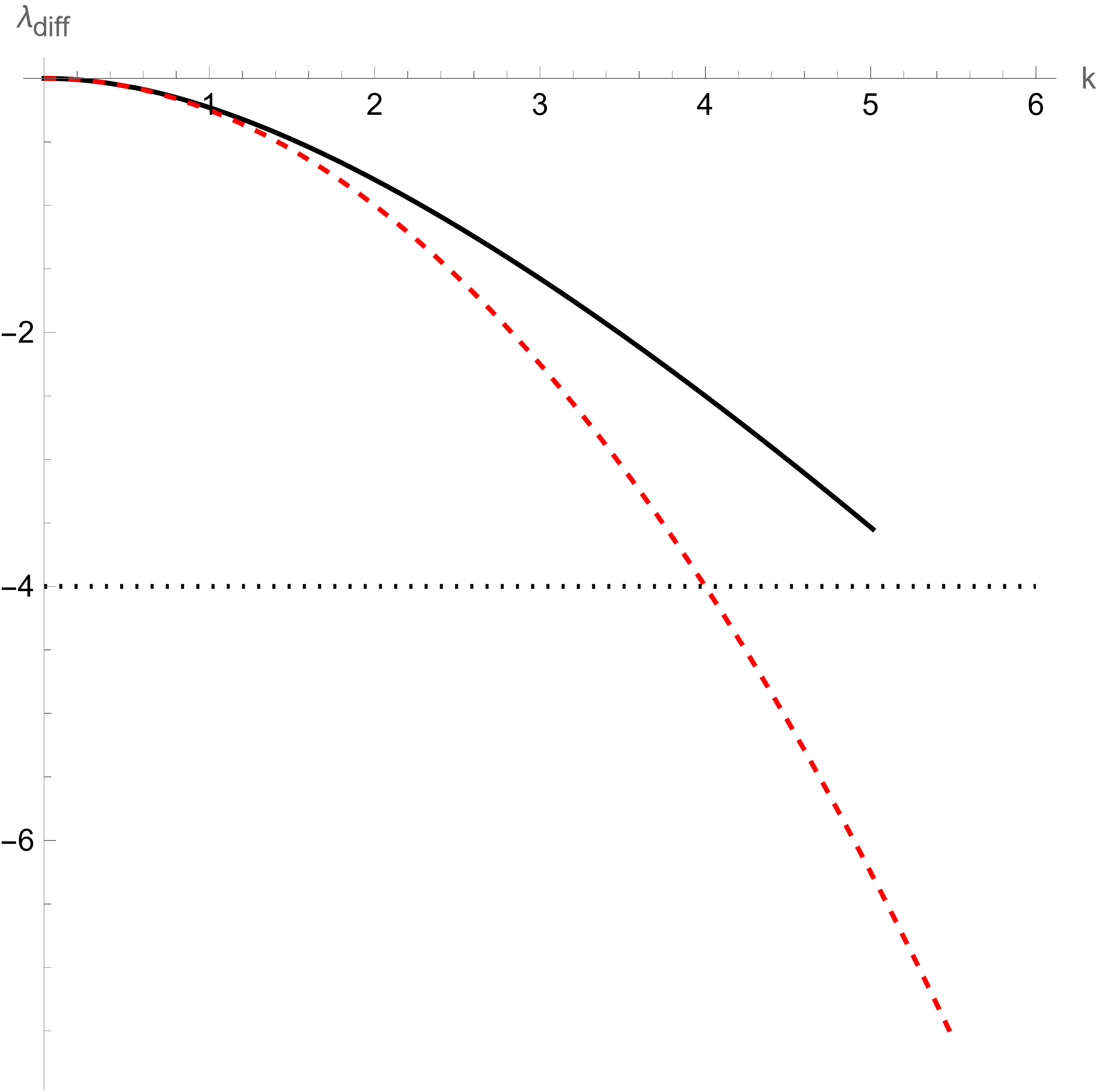}
		\caption{Diffusion mode.}
	\end{subfigure}%
	\begin{subfigure}{.5\textwidth}
		\centering
		\includegraphics[width=0.8\linewidth]{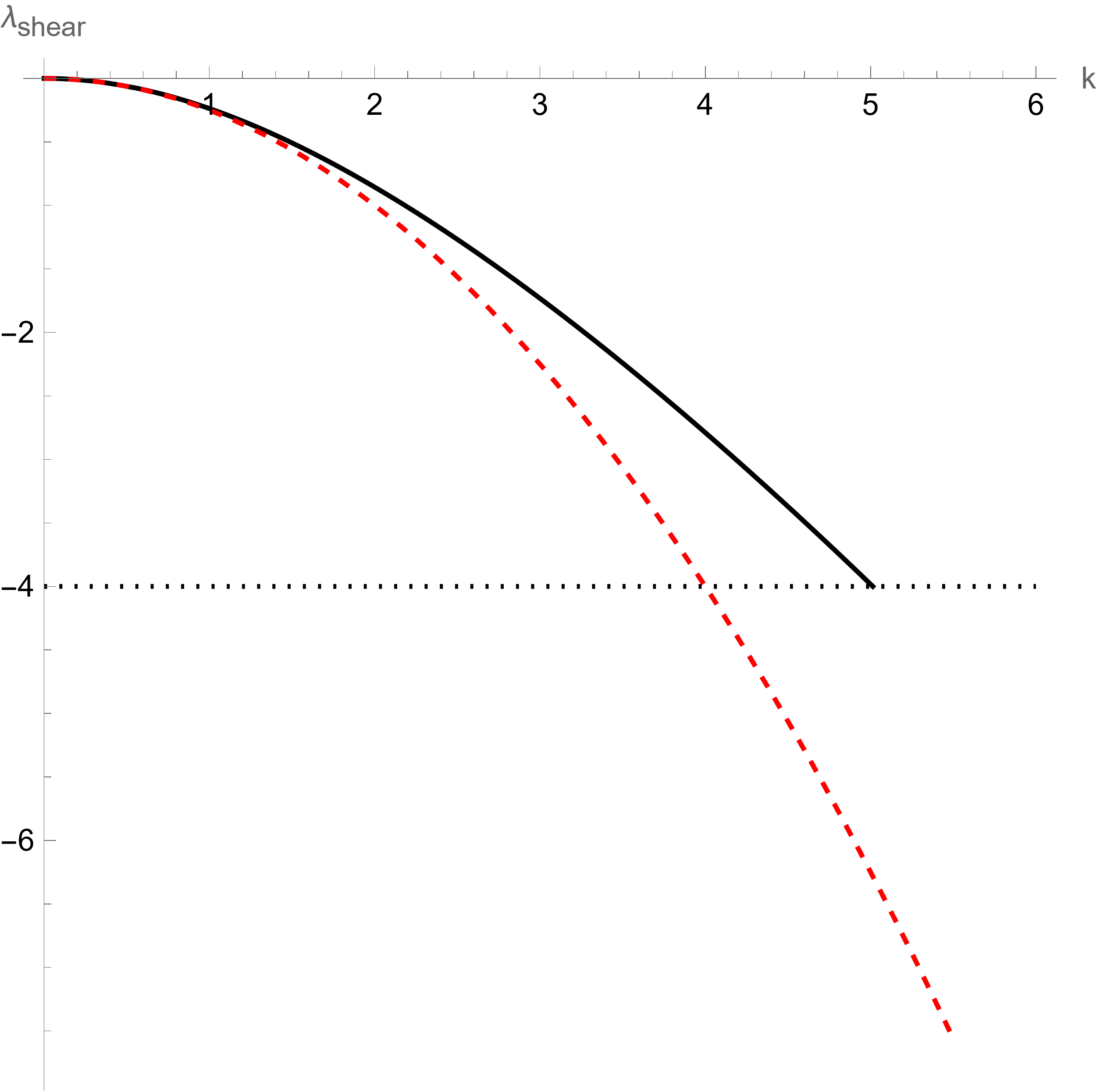}
		\caption{Shear mode.}
	\end{subfigure}
 \begin{subfigure}{.5\textwidth}
		\centering
		\includegraphics[width=0.8\linewidth]{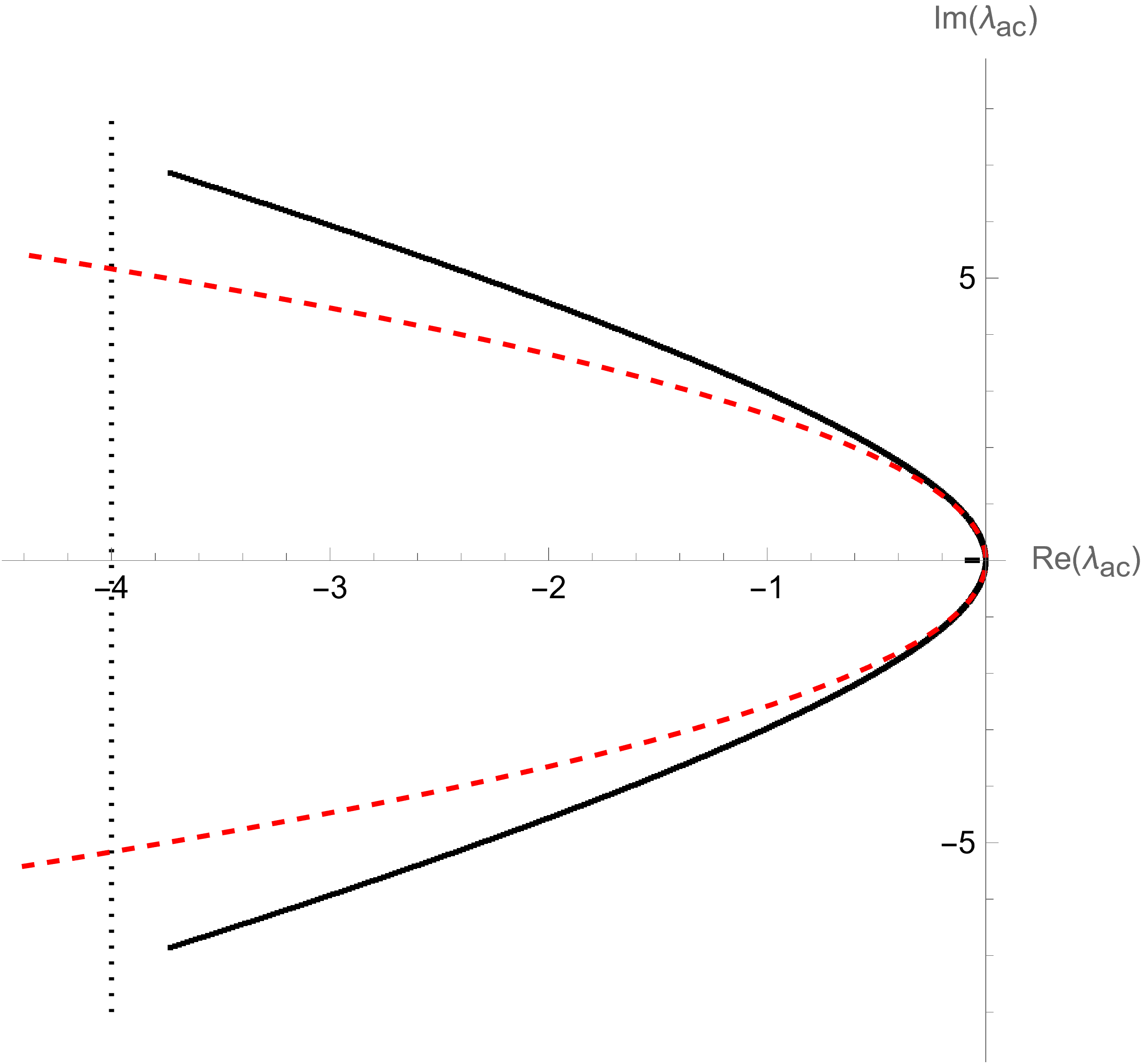}
		\caption{Acoustic mode.}
	\end{subfigure}
    \caption{Modal branches in dependence on the wave number (black solid line) compared to the leading-order polynomial approximation (red dashed line) for  $\tau=0.25$ and the essential spectrum at $-\frac{1}{\tau}$ (dotted black line). The hydrodynamic branches terminate at the minimal critical wave number $ k_{\rm crit, min}$ \eqref{kcritmin}.}
    \label{Modes}
\end{figure}

%
%
Critical wave numbers for the branches are found to be \cite{kogelbauer2023exact},
\begin{equation}\label{kcrit}
    \begin{split}
       & k_{\rm crit}(\lambda_{\rm shear})  = \sqrt{\frac{\pi}{2}}\frac{1}{\tau}\approx 1.2533\frac{1}{\tau},\\
       &  k_{\rm crit}(\lambda_{\rm diff})  \approx 1.3560 \frac{1}{\tau},\\
       &  k_{\rm crit}(\lambda_{\rm ac}){= k_{\rm crit}(\lambda^*_{\rm ac})}  \approx 1.3118 \frac{1}{\tau},
    \end{split}
\end{equation}
In particular, we note that the critical wave number of the shear mode is minimal, which implies that all three branches exists for $0<k<k_{\rm crit}$ and we set
\begin{equation}\label{kcritmin}
    k_{\rm crit, min} = k_{\rm crit}(\lambda_{\rm shear}) = \sqrt{\frac{\pi}{2}}\frac{1}{\tau}.
\end{equation}

The eigenvalues admit the following asymptotic expansions in terms of the wave number:
\begin{equation}\label{expandlambda}
\begin{split}
    \lambda_{\rm diff}(k) & = -\tau k^2+\frac{9}{5}\tau^3k^4+\mathcal{O}(k^6),\\
    \lambda_{\rm shear}(k) & = -\tau k^2+\tau^3 k^4+\mathcal{O}(k^6),\\
    \lambda_{\rm ac}(k) & = \ri\sqrt{\frac{5}{3}}k -\tau k^2 +\ri\frac{7\tau^2}{6\sqrt{15}}k^3+\frac{62}{45}\tau^3k^4+\mathcal{O}(k^5),
    \end{split}
\end{equation}
which can be seen from Taylor expanding $\lambda$ in $k$ and comparing powers in \eqref{defSigma}, see also \cite{kogelbauer2023exact}. 
\begin{remark}
    Since zero is a five-fold degenerate eigenvalue for $k=0$, we do not expect - in general - that the eigenvalues depend analytically on $k$. Spectral perturbation theory only guarantees the expansion in a Puiseux series, i.e., a Taylor series in $k^{{1}/{5}}$. For the BGK equation, however, the fractional terms cancel out and only powers in $k$ remain, which is consistent with Ellis \& Pinsky  \cite{ellis1975first}. 
\end{remark}
\begin{remark}
	While the structure of the Boltzmann-BGK spectrum \eqref{thmspec} agrees with and is a special case of a more general linear Boltzmann equation \cite{dudynski2013,palczewski1984}, the knowledge of the spectral function \eqref{defSigma} allows to discern more detailed analytical information about the hydrodynamic spectrum, in particular, the accurate estimate for the critical wave numbers \eqref{kcrit}.
\end{remark}

 Figure \ref{Modes} shows the dependence of the modes on wave number in comparison to their leading-order polynomial approximation in \eqref{expandlambda}, which correspond to Euler and Navier--Stokes equations.

 There exists exactly five discrete, isolated eigenvalues of the operator $\mathcal{L}_{\mathbf{k}}$ with an associated five-dimensional eigenspace, which we will call the \textit{hydrodynamic manifold}. This manifold will serve as our constraint to define a closure operator. Because the five eigenvalues are above the essential spectrum, any solution restricted to the given wave number, will decay exponentially fast to the hydrodynamic manifold, rendering it a \textit{slow manifold}.  

In the following, we will be interested in the linear subspace generated by the eigenfunctions associated to $\Lambda_{\rm BGK} = \{\lambda_{\rm diff},\lambda_{\rm ac},\lambda_{\rm ac}^{*}, \lambda_{\rm shear}\}$ at each wave number $k$. To ease notation, we bundle the eigenvalues in a vector (counted with multiplicity):
\begin{equation}
    \boldsymbol{\lambda} = (\lambda_{\rm diff},\lambda_{\rm ac},\lambda_{\rm ac}^*,\lambda_{\rm shear},\lambda_{\rm shear}),
\end{equation}
and define the diagonal matrix
\begin{equation}
    \boldsymbol{\Lambda} = \text{diag}(\boldsymbol{\lambda}). 
\end{equation}
Also, we denote the set and vector of simple eigenvalues as
\begin{equation}
    \Lambda_{\rm simple} = \{\lambda_{\rm diff},\lambda_{\rm ac},\lambda_{\rm ac}^*\},\quad \boldsymbol{\lambda}_{\rm simple} = (\lambda_{\rm diff},\lambda_{\rm ac},\lambda_{\rm ac}^*).
\end{equation}

For each wave vector $\mathbf{k}$ with $0<k<k_{\rm crit,min}$, the eigenspace associated to the modes spans a five-dimensional linear subspace, which we call the \textit{hydrodynamic manifold}:   
\begin{equation}\label{defhydro}
    \mathcal{M}_{\rm hydro,\mathbf{k}} = \text{span}_{\lambda\in\Lambda_{\rm simple}} \hat{f}_{\lambda}(\mathbf{k}) \otimes \text{span} \{\hat{f}_{\lambda_{\rm shear},1}(\mathbf{k}),\hat{f}_{\lambda_{\rm shear},2}(\mathbf{k}) \}. 
\end{equation}
The hydrodynamic manifold \eqref{defhydro} is invariant with respect to the flow generated by \eqref{maineq}. We write
\begin{equation}\label{dynfhydro}
\hat{f}_{\rm hydro}(\mathbf{k},\mathbf{v},t) = \sum_{\lambda\in \Lambda_{\rm simple}} \alpha_{\lambda}(t)  {\hat{f}}_{\lambda}(\mathbf{v},\mathbf{k}) 
+ \alpha_{\lambda_{\rm shear},1}(t)\hat{f}_{\lambda_{\rm shear},2}(\mathbf{v},\mathbf{k})
+\alpha_{\lambda_{\rm shear},2}(t)\hat{f}_{\lambda_{\rm shear},2}(\mathbf{v},\mathbf{k}),
\end{equation}
for a solution on $\mathcal{M}_{\rm hydro }$. The vector
\begin{equation}
\boldsymbol{\alpha} = (\alpha_{\lambda_{\rm diff}},{\alpha_{\lambda_{\rm ac}},\alpha_{\lambda^*_{\rm ac}},\alpha_{\lambda_{\rm shear, 1}},\alpha_{\lambda_{\rm shear,2}}}),
\end{equation}
is comprised of \textit{spectral variables} or \textit{spectral coefficients}.

\subsection{Spectral Closure for the Boltzmann-BGK Equation}
Given the hydrodynamic manifold $\mathcal{M}_{\rm hydro,\mathbf{k}}$ for $0<k<k_{\rm crit,min}$, which is spanned by the eigenvectors $\hat{f}_{\lambda}$, we define the \textit{spectral closure} as, 
\begin{equation}\label{defCspectral}
\begin{split}
    \mathcal{C}_{\rm spectral}&: \text{range } \mathbb{P}_5|_{\mathcal{M}_{\rm hydro,\mathbf{k}}}\to \text{range } \mathbb{P}_5^{\perp }|_{\mathcal{M}_{\rm hydro,\mathbf{k}}},\\
    \mathcal{C}_{\rm spectral}&(\mathbb{P}_5\hat{f}_\lambda) = \mathbb{P}_5^{\perp}\hat{f}_{\lambda}. 
  \end{split}  
\end{equation}
The closure operator \eqref{defCspectral}, defined only on the $\mathcal{M}_{\rm hydro,\mathbf{k}}$, maps - for each eigenvalue - the first five moments of an eigenvector to the orthogonal complement on the same eigenvector. The closure operator \eqref{defCspectral} is defined with respect to the spectral basis, whereas the closure formalism \eqref{oneplusC} is defined with respect to macroscopic variables. To obtain the corresponding change of coordinates, let us denote the first five elements of a simple eigenfunction $\hat{f}_{\lambda}$ as
\begin{equation}
    \boldsymbol{\eta}(\lambda) = \mathbb{P}_5\hat{f}_{\lambda},
\end{equation}
while we write 
\begin{equation}
    \boldsymbol{\eta}_1(\lambda_{\rm shear}) = \mathbb{P}_5\hat{f}_{\lambda_{shear},1},\qquad \boldsymbol{\eta}_2(\lambda_{\rm shear}) = \mathbb{P}_5\hat{f}_{\lambda_{shear},2}.
\end{equation}
Taking projections in \eqref{dynfhydro}, we have 
\begin{equation}\label{defhhydro}
    \hat{\mathbf{h}}_{\rm hydro}(\mathbf{v},t) = \sum_{\lambda\in \Lambda_{\rm simple}} \alpha_{\lambda}(t) \mathbb{P}_5 \hat{f}_{\lambda}(\mathbf{v},\mathbf{k}) + \alpha_{\lambda_{\rm shear},1}\mathbb{P}_5\hat{f}_{\lambda_{\rm shear},1}(\mathbf{v},\mathbf{k})+\alpha_{\lambda_{\rm shear},2}\mathbb{P}_5\hat{f}_{\lambda_{\rm shear},2}(\mathbf{v},\mathbf{k}),
\end{equation}
where we have suppressed the explicit dependence on the wave vector $\mathbf{k}$. 
To ease notation, we define the $5\times 5$ matrix of spectral basis vectors for the BGK equation (see Theorem \ref{thmspec}) as
\begin{equation}\label{H}
    \mathbf{H} := [\boldsymbol{\eta}(\lambda_{\rm diff}),\boldsymbol{\eta}(\lambda_{\rm ac}),\boldsymbol{\eta}(\lambda_{ac}^*),\boldsymbol{\eta}_1(\lambda_{\rm shear}),\boldsymbol{\eta}_2(\lambda_{\rm shear}) ],
\end{equation}
which allows us to write the macroscopic variables on the hydrodynamic manifold as
\begin{equation}\label{defhhydroBGK}
    \hat{\mathbf{h}}_{\rm hydro} =\mathbf{H}\boldsymbol{\alpha}.
\end{equation}
Since \eqref{defhhydroBGK} is composed of eigenvectors entirely, the evolution on the hydrodynamic manifold in terms of spectral variables simply becomes
\begin{equation}\label{dynalpha}
  \frac{d \boldsymbol{\alpha}}{dt} = \boldsymbol{\Lambda}\boldsymbol{\alpha},  
\end{equation}
i.e., the spectral dynamics diagonalize completely since geometric and algebraic multiplicity are equal for each mode.\\

To define the spectral closure for macroscopic variables, we define
\begin{equation}\label{Flambda}
    F_{\boldsymbol{\lambda}} = [\hat{f}_{\lambda_{\rm diff}},\hat{f}_{\lambda_{\rm ac}},\hat{f}_{\lambda_{\rm ac}^*},\hat{f}_{\lambda_{\rm shear},1},\hat{f}_{\lambda_{\rm shear},2}],
\end{equation}
which implies that $\mathbf{H}=\mathbb{P}_5F_{\boldsymbol{\lambda}}$. Based on the coordinate change \eqref{H} and \eqref{Flambda}, the closure operator in macroscopic variables then reads
\begin{equation}\label{closureBGK}
\begin{split}
    \mathcal{C}_{\rm spectral}\hat{\mathbf{h}}_{\rm hydro} & = \mathcal{C}_{\rm spectral}\mathbf{H}\boldsymbol{\alpha}\\
    & = \mathbb{P}_5^{\perp}F_{\boldsymbol{\lambda}}\boldsymbol{\alpha}\\
    & = \mathbb{P}_5^{\perp}F_{\boldsymbol{\lambda}}\mathbf{H}^{-1}\hat{\mathbf{h}}_{\rm hydro}
    \end{split}
\end{equation}
provided that the inverse exists (which will be elaborated in Section \ref{HydroEqu}). \\
We emphasize that the spectral closure for the BGK equation \eqref{closureBGK} is an exact (invariant) closure in the sense of \eqref{exactclosurek} by construction. Indeed, we find that evaluating \eqref{exactclosurek} on $\hat{\mathbf{h}}_{\rm hydro}$ gives:
\begin{equation}
\begin{split}
    (\mathcal{C}_{\rm spectral}\mathbb{P}_5-\mathbb{P}_5^{\perp})\mathcal{L}_{\mathbf{k}}(1+\mathcal{C}_{\rm spectral})\hat{\mathbf{h}}_{\rm hydro} & =  (\mathcal{C}_{\rm spectral}\mathbb{P}_5-\mathbb{P}_5^{\perp})\mathcal{L}_{\mathbf{k}}(1+\mathcal{C}_{\rm spectral})\mathbf{H}\boldsymbol{\alpha}\\
    & = (\mathcal{C}_{\rm spectral}\mathbb{P}_5-\mathbb{P}_5^{\perp})\mathcal{L}_{\mathbf{k}}F_{\boldsymbol{\lambda}}\boldsymbol{\alpha}\\
    & = (\mathcal{C}_{\rm spectral}\mathbb{P}_5-\mathbb{P}_5^{\perp})F_{\boldsymbol{\lambda}}\boldsymbol{\Lambda}\boldsymbol{\alpha}\\
    & = (\mathcal{C}_{\rm spectral}\mathbf{H}-\mathbb{P}_5^{\perp}F_{\boldsymbol{\lambda}})\boldsymbol{\Lambda}\boldsymbol{\alpha}\\
    & =0,
    \end{split}
\end{equation}
where in the third step, we have used that columns of $F_{\boldsymbol{\lambda}}$ are eigenvectors of $\mathcal{L}_{\mathbf{k}}$. 


\begin{remark}
Because of the existence of a critical wave number for each branch of eigenvalues \eqref{kcrit}, a full set of five eigenvalues only exists up to $k_{\rm crit,min}$. For $k>k_{\rm crit,min}$, the modes vanish one by one, which implies that the full set of five macroscopic variables cannot be resolved uniquely any longer. In particular, the matrix $\mathbf{H}$ might not be defined as a square matrix \eqref{H}, but much rather as a rectangular matrix. Also, the inverse appearing in \eqref{closureBGK} has to be understood as a generalized inverse (e.g., pseudo-inverse). Implications of this degeneracy shall not be further discussed in this paper.
\end{remark}

Next section will be devoted to the explicit calculation of the basis matrix $\mathbf{H}$ (the invertability of $\mathbf{H}$ will be discussed in Section \ref{HydroEqu}). Through a simple linear change of coordinates, we then obtain the dynamics for the macroscopic moments.

\section{From Spectral Coordinates to Macroscopic Variables}
\label{sec:spectral_closure}
In this section, we construct the exact spectral closure for the BGK equation based on the knowledge of the spectrum \eqref{thmspec}. We derive the coordinate change from spectral parameters to the primitive, macroscopic variables (i.e., the basis matrix $\mathbf{H}$) in two consistency ways: First, we derive a general algebraic from of the first five moments of a simple eigenfunction from the interplay of the linear transport and collision operators. This analysis will be specific for the BGK equation and depends on the specific form of the projection operator \eqref{defP5}.\\
Secondly, we use analytical spectral calculus and Riesz projections to obtain the same result and to show consistency of the two approaches. We emphasize that the approach via spectral projections, although equivalent, can be be applied to a more general setting as well.\\
Before we proceed, let us collect some notation and results from \cite{kogelbauer2023exact} regarding the spectral problem of the linearized Boltzmann-BGK operator \eqref{defLhat}, which will be useful in the following two subsections.
We define the Green's function matrices as
\begin{equation}\label{defGST}
\begin{split}
G_{{{L}}}(z,n,m) &= \langle (\ri \tau \mathbf{v}\cdot\mathbf{k}-\mathbb{P}_5-z)^{-1}e_n,e_m\rangle_{\mathbf{v}},\\
G_{{S}}(z,n,m) &= \langle (\ri\tau\mathbf{v}\cdot\mathbf{k}-z)^{-1}e_n,e_m\rangle_{\mathbf{v}},
\end{split}
\end{equation}
which satisfy the equation \cite{kogelbauer2023exact},
\begin{equation}\label{Greenseq}
G_{{{L}}}=G_{{S}}+G_{{{L}}}G_{{S}}.
\end{equation}
 Assuming $\det(\Id -G_{{S}})\neq 0$, equation \eqref{Greenseq} can be solved to get
\begin{equation}\label{GTGS}
G_{{{L}}}=G_{{S}}(\Id -G_{{S}})^{-1} = (\Id -G_{{S}})^{-1} - \Id. 
\end{equation}
From \cite{kogelbauer2023exact} we know that the matrix $z\mapsto G_{{S}}(z)$ can be conjugated using the rotation matrix \eqref{defQtilde}: 
\begin{equation}\label{GsmI}
    \left. G_{{S}}(z)\right|_{z=\ri k\tau \zeta}-\Id = \frac{1}{\ri \tau k}\tilde{\mathbf{Q}}_{\mathbf{k}} (G(\zeta)-\ri \tau k)\tilde{\mathbf{Q}}_{\mathbf{k}}^T,
\end{equation}
where the matrix $G(\zeta)$ reads,
\begin{equation}\label{defGzeta}
G(\zeta) = \begin{pmatrix}
Z(\zeta) & 1+\zeta Z(\zeta) & 0 & 0 & \frac{\zeta +(\zeta^2-1)Z(\zeta)}{\sqrt{6}}\\
1+\zeta Z(\zeta) & \zeta + \zeta^2 Z(\zeta) & 0 & 0 &   \frac{\zeta^2+(\zeta^3-\zeta)Z(\zeta)}{\sqrt{6}}\\
0 & 0 &  Z(\zeta) & 0 & 0\\
0 & 0 &  0 & Z(\zeta) & 0\\
\frac{\zeta +(\zeta^2-1)Z(\zeta)}{\sqrt{6}} & \frac{\zeta^2+(\zeta^3-\zeta)Z(\zeta)}{\sqrt{6}}& 0 & 0 & \frac{\zeta^3-\zeta+(\zeta^4-2\zeta^2+5)Z(\zeta)}{6}
\end{pmatrix},
\end{equation}
while $\zeta\mapsto Z(\zeta)$ is the plasma dispersion function \eqref{defZ}. 
Furthermore, the spectral function \eqref{defSigma} is related to Green's matrix $G_{{S}}$ \eqref{defGST} as \cite{kogelbauer2023exact},
\begin{equation}\label{SigmaG}
    \Sigma_{k,\tau}(\zeta) = \left. \det(G_{{S}}(z)-\Id)\right|_{z=\ri k\tau \zeta}. 
\end{equation}
With  \eqref{GsmI}, \eqref{defGzeta}, \eqref{defQtilde} and \eqref{defQ}, the determinant in \eqref{SigmaG} is evaluated easily to get the closed-form expression for the spectral function \eqref{defSigma} obtained in  \cite{kogelbauer2023exact}.

\subsection{Spectral-to-Hydrodynamic Coordinate Transform by Spectral Temperature}
\label{sec:spectral_temperature}

An eigenvector $\hat{f}_{\lambda}$ of \eqref{defLhat} with eigenvalue $\lambda$ satisfies the equation
\begin{equation}\label{eigf1}
    -\ri\mathbf{v}\cdot\mathbf{k} \hat{f}_{\lambda} -\frac{1}{\tau}\hat{f}_{\lambda} +\frac{1}{\tau}\mathbb{P}_5 \hat{f}_{\lambda} = \lambda \hat{f}_{\lambda},
\end{equation}
or, equivalently,
\begin{equation}\label{flambdaeig}
    \hat{f}_{\lambda} = \frac{\mathbb{P}_5\hat{f}_{\lambda}}{\tau\lambda + 1 + \ri\tau\mathbf{v}\cdot\mathbf{k}}.
\end{equation}
We emphasize that the numerator in \eqref{flambdaeig} is always non-zero for the range of $k$ for which $\lambda$ is defined since
\begin{equation}
-\frac{1}{\tau}<\Re\lambda(k) <0,
\end{equation}
for $0<k<k_{\rm crit,min}$.
Projecting equation \eqref{flambdaeig} via $\mathbb{P}_5$ gives the following implicit equation for the first five entries of an eigenvector:
\begin{equation}\label{etaeq}
    \boldsymbol{\eta}(\lambda) = \left\langle \frac{\mathbf{e}(\mathbf{v})\cdot \boldsymbol{\eta}(\lambda)}{\tau\lambda + 1 + \ri\tau\mathbf{v}\cdot\mathbf{k}},\mathbf{e}(\mathbf{v})\right\rangle_{\mathbf{v}}. 
\end{equation}
Writing $\boldsymbol{\eta} = (\eta_1,\eta_2,\eta_3,\eta_4,\eta_5)$ and integrating \eqref{eigf1} over $\mathbf{v}$ implies the relation
\begin{equation}
   -\ri \mathbf{k}\cdot(\eta_1,\eta_2,\eta_3) =  \lambda \eta_1,
\end{equation}
or, equivalently, in terms of the splitting \eqref{split}: 
\begin{equation}\label{formetas}
    (\eta_1,\eta_2,\eta_3) = \frac{\ri\lambda \eta_1}{k^2}\mathbf{k} + (\eta_1,\eta_2,\eta_3)_{\perp}. 
\end{equation}
Using the Green's functions matrices \eqref{defGST}, we can rewrite \eqref{etaeq} as a functional eigenvalue problem:
\begin{equation}
     \boldsymbol{\eta}(\lambda) = G_{{S}}(-\tau\lambda-1) \boldsymbol{\eta}(\lambda),
\end{equation}
which, in the light of \eqref{GsmI}, can be rewritten as
\begin{equation}
    \tilde{\mathbf{Q}}_{\mathbf{k}}^T\boldsymbol{\eta} \in \ker(G(\zeta)-\ri k \tau),
\end{equation}
for $\ri(1+\tau\lambda)=k \tau \zeta$. Given the structure of \eqref{defGzeta} and the definition of the modes through \eqref{defSigma}, we immediately see that
\begin{equation}\label{basisshear}
\boldsymbol{\eta}_1(\lambda_{shear})=\tilde{\mathbf{Q}}_{\mathbf{k}}\left(\begin{array}{c}
0\\
0\\
1\\
0\\
0
\end{array}
\right),\qquad \boldsymbol{\eta}_2(\lambda_{shear})=\tilde{\mathbf{Q}}_{\mathbf{k}}\left(\begin{array}{c}
0\\
0\\
0\\
1\\
0
\end{array}
\right),
\end{equation}
since the middle block in $G(\zeta)$ decouples from the other part of the matrix (which of course corresponds to the factorisation in \eqref{defSigma}).\\
To obtain the structure of the simple eigenvectors, we first assume, by rescaling $\hat{f}_{\lambda}$ accordingly, that
\begin{equation}
    \langle \hat{f}_{\lambda}, 1 \rangle_{\mathbf{v}} = 1.
\end{equation}
We note that, using \eqref{formetas}, that columns one two and five of \eqref{defGzeta} are scalar multiples of each other for $(\eta_1,\eta_2,\eta_3)_{\perp}=0$. Consequently, for a simple eigenvalue $\lambda$, we can set 
\begin{equation}\label{hkexplicit}
    \boldsymbol{\eta}(\lambda) = \left(\begin{array}{c} 1 \\ \frac{\ri\lambda}{k^2}\mathbf{k}\\ \theta(\lambda)\end{array}\right),
\end{equation}
for some function $\lambda\mapsto \theta(\lambda)$. We call the basis vectors \eqref{hkexplicit} for a simple eigenvalue \textit{k-aligned}.\\
Since the two eigenvectors for the shear mode are complete explicit \eqref{basisshear}, all non-trivial information of the spectral closure is encoded in the function $\lambda\mapsto\theta(\lambda)$, which we call \textit{spectral temperature}. Since we know the structure of a simple eigenvector \eqref{hkexplicit}, we can derive a formula for the spectral temperature as follows. Using \eqref{hkexplicit}, we take an inner product of \eqref{flambdaeig} with $e_4$ to find
\begin{equation}
    \theta(\lambda) = \left\langle\frac{1+\ri\frac{\lambda}{k^2}\mathbf{k}\cdot\mathbf{v}+\theta(\lambda)e_4(\mathbf{v})}{\tau\lambda + 1 + \ri\tau\mathbf{v}\cdot\mathbf{k}}, e_4(\mathbf{v})\right\rangle_{\mathbf{v}},
\end{equation}
which can be solved to 
\begin{equation}\label{theta1}
    \theta(\lambda) = \frac{\frac{1}{k^2}\left\langle\frac{k^2+\ri\lambda\mathbf{k}\cdot\mathbf{v}}{\tau\lambda + 1 + \ri\tau\mathbf{v}\cdot\mathbf{k}},e_4(\mathbf{v})\right\rangle_{\mathbf{v}}}{1-\left\langle\frac{e_4(\mathbf{v})}{\tau\lambda + 1 + \ri\tau\mathbf{v}\cdot\mathbf{k}},e_4(\mathbf{v})\right\rangle_{\mathbf{v}}}. 
\end{equation}
Expression \eqref{theta1} can then be evaluated explicitly for each simple $\lambda$. We refer to Subsection \ref{spectraltemperature} for an explicit formula and properties of $\theta$.\\ 
In the next subsection, we will give an alternative derivation of \eqref{hkexplicit} and \eqref{basisshear} using spectral projections to emphasize consistency.

\subsection{Spectral-to-Hydrodynamic Coordinate Transform by Riesz Projections}
\label{sec:Riesz}
In the following, we show that the the spectral basis \eqref{basisshear} and \eqref{hkexplicit} obtained in the previous section can be derived equivalently through spectral calculus. Indeed, for any set of discrete, isolated eigenvalues $k\mapsto\Lambda_{\rm BGK}(k)\subset\mathbb{C}$, depending on wave number, we can define the \textit{Riesz projection} as 
\begin{equation}\label{defPLambda}
     \mathbb{P}_{\Lambda} = -\frac{1}{2\ri\pi}\oint_{\Gamma(\Lambda_{\rm BGK})}(\mathcal{L}_{\mathbf{k}}-w)^{-1}\, dw,
\end{equation}
where $\Gamma(\Lambda_{\rm BGK})$ is a simple contour in the complex plain, encircling the full spectral set $\Lambda_{\rm BGK}=\{\lambda_{\rm diff},\lambda_{\rm ac},\lambda_{\rm ac}^*,\lambda_{\rm shear}\}$ once in positive direction. From analytical spectral calculus \cite{hislop2012introduction}, we know that \eqref{defPLambda} is indeed a projection, whose range is given by the invariant subspace (generalized eigenspace) associated to $\Lambda_{\rm BGK}$. In particular, we see from \eqref{defhhydro}, it follows that
\begin{equation}\label{defH}
    \mathbf{H} = -\frac{1}{2\ri\pi}\mathbb{P}_5\oint_{\Gamma(\Lambda_{\rm BGK})} (\mathcal{L}_{\mathbf{k}}-w)^{-1}\, dw\, \mathbb{P}_5.
\end{equation}
\begin{remark}
Here, we have assumed that the five basis vectors $\mathbf{e}$ are mapped to five linearly independent vectors on $\mathcal{M}_{\rm hydro}(\Lambda)$ (which is indeed the case for the BGK equation). For a more general kinetic model, this might not be the case and a non-invertability of the spectral basis in terms of the macroscopic variables would indicate a restriction of the hydrodynamic dynamics for the given range of wave numbers. 
\end{remark}
Setting $w=-\frac{1}{\tau}(z+1)$, the resolvent transforms according to
\begin{equation}\label{reszw}
\begin{split}
(\mathcal{L}_{\mathbf{k}}-w)^{-1} & = \left(-\ri\mathbf{v}\cdot\mathbf{k}-\frac{1}{\tau}+\frac{1}{\tau}\mathbb{P}_5+\frac{1}{\tau}(z+1)\right)^{-1}\\
& =-\tau \left(\ri\tau\mathbf{v}\cdot\mathbf{k}-\mathbb{P}_5-z\right)^{-1}.
\end{split}
\end{equation}
Using the second resolvent identity together with \eqref{reszw} and $dw=-\frac{1}{
\tau}dz$, we can then write:
\begin{equation}
    \begin{split}
   \mathbb{P}_{\Lambda}e_m  &= -\frac{1}{2\ri\pi}\oint_{\Gamma(\Lambda_{\rm BGK})}(\mathcal{L}_{\mathbf{k}}-w)^{-1}e_m\, dw\\
   & = -\frac{1}{2\ri\pi}\oint_{\Gamma(\Lambda_\tau)}(\ri\tau\mathbf{v}\cdot\mathbf{k}-\mathbb{P}_5-z)^{-1}e_m\, dz\\
   & =   -\frac{1}{2\ri\pi}\oint_{\Gamma(\Lambda_\tau)}(\ri\tau\mathbf{v}\cdot\mathbf{k}-z)^{-1}e_m+ (\ri\tau\mathbf{v}\cdot\mathbf{k}-z)^{-1}\mathbb{P}_5(\ri\tau\mathbf{v}\cdot\mathbf{k} - \mathbb{P}_5 - z)e_m\, dz\\
   & =  -\frac{1}{2\ri\pi}\oint_{\Gamma(\Lambda_\tau)}(\ri\tau\mathbf{v}\cdot\mathbf{k}-z)^{-1}\sum_{j=0}^4\langle(\ri\tau\mathbf{v}\cdot\mathbf{k}-\mathbb{P}_5-z)^{-1}e_m,e_j\rangle_{\mathbf{v}}e_j\, dz\\
&=-\frac{1}{2\ri\pi}\oint_{\Gamma(\Lambda_\tau)}\sum_{j=0}^4G_{{L}}(z,m,j)(\ri\tau\mathbf{v}\cdot\mathbf{k}-z)^{-1}e_j dz,
    \end{split}
\end{equation}
where we have set
\begin{equation}
\Gamma(\Lambda_{\rm BGK}) = -\frac{1}{\tau} (\Gamma(\Lambda_{\tau})+1).
\end{equation}
With the notation \eqref{defR}, we have that $w\in\mathcal{R}_a \iff z \in \mathcal{R}_{1}$. Using relation \eqref{Greenseq} between the Green's function matrices $G_{{S}}$ and $G_{{L}}$ together with the fact that $z\mapsto G_{{S}}(z)$ is holomorphic in $\mathcal{R}_1$, we arrive at 
\begin{equation}
    \begin{split}
    \langle \mathbb{P}_\lambda e_m,e_n\rangle_{\mathbf{v}} & =-\frac{1}{2\ri\pi}\oint_{\Gamma(\Lambda_\tau)}\sum_{j=0}^4G_{{L}}(z,m,j)\langle(\ri\tau\mathbf{v}\cdot\mathbf{k}-z)^{-1}e_j,e_n\rangle_{\mathbf{v}} dz\\
    & = -\frac{1}{2\ri\pi}\oint_{\Gamma(\Lambda_\tau)}\sum_{j=0}^4G_{{L}}(z,m,j)G_{{S}}(z,j,n) dz\\
    &=  -\frac{1}{2\ri\pi}\oint_{\Gamma(\Lambda_\tau)} (G_{{L}}(z)G_{{S}}(z))_{n,m} dz\\
    &= -\frac{1}{2\ri\pi}\oint_{\Gamma(\Lambda_\tau)} G_{{L}}(z,n,m) dz\\
    &= \frac{1}{2\ri\pi}\oint_{\Gamma(\Lambda_\tau)} [(G_{{S}}(z)-\Id)^{-1}]_{n,m} dz.
    \end{split}
\end{equation}
By relation \eqref{SigmaG} and by applying the Residue Theorem, we find that 
\begin{equation}\label{res1}
\begin{split}
    \langle \mathbb{P}_\lambda e_m,e_n\rangle_{\mathbf{v}} & = -\frac{1}{2\ri\pi}\oint_{\Gamma(\Lambda_\tau)} [(\Id-G_{{S}}(z))^{-1}]_{n,m} dz\\
    &=-\frac{1}{2\ri\pi}\oint_{\Gamma(\Lambda_\tau)} \Sigma_{k,\tau}^{-1}(z)\adj(\Id - G_{{S}}(z)) dz\\
    & = \sum_{\lambda_{\tau}\in\Lambda_\tau}\mathrm{Res}_{z\to\lambda_\tau} \Sigma_{k,\tau}^{-1}(z)\adj(G_{{S}}(z)-\Id).
\end{split}
\end{equation}
For a simple eigenvalue $\lambda_{\tau}$, the function $\Sigma_{k,\tau}^{-1}$ has a pole of order one at $\lambda_{\tau}$ and 
\begin{equation}
    \langle \mathbb{P}_\lambda e_m,e_n\rangle_{\mathbf{v}} \sim \adj( G_{{S}}(\lambda_{\tau})-\Id),
\end{equation}
where $\sim$ indicates equality up to multiplication by a complex number, while at the shifted shear mode $\lambda_{\text{shear},\tau}$, the function $\Sigma_{k,\tau}^{-1}$ has a pole of order two and 
\begin{equation}
\begin{split}
    \langle \mathbb{P}_\lambda e_m,e_n\rangle_{\mathbf{v}} & = \mathrm{Res}_{z\to\lambda_{\text{shear},\tau}} \Sigma_{k,\tau}^{-1}(z)\adj(G_{{S}}(z)-\Id)\\
     & = \lim_{z\to\lambda_{\text{shear},\tau}}\frac{d}{dz}\Big[(z-\lambda_{\text{shear},\tau})^2\Sigma_{k,\tau}^{-1}(z)\adj(G_{{S}}(z)-\Id)\Big]\\
     &=\lim_{z\to\lambda_{\text{shear},\tau}}\Big[(z-\lambda_{\text{shear},\tau})^2\Sigma_{k,\tau}^{-1}(z)\Big]'\adj(G_{{S}}(\lambda_{\text{shear},\tau})-\Id)\\
     &\qquad+\lim_{z\to\lambda_{\text{shear},\tau}}(z-\lambda_{\text{shear},\tau})^2\Sigma_{k,\tau}^{-1}(z)\Big[\adj( G_{{S}}(z)-\Id)\Big]'\\
     & \sim \lim_{z\to\lambda_{\text{shear},\tau}} \Big[\adj( G_{{S}}(z)-\Id)\Big]',
    \end{split}
\end{equation}
where prime denotes the derivative $d/dz$, and where we have used that $\adj (A)=0$ if $\dim\ker A\geq 2$. 
The formula 
\begin{equation}
\adj(AB)=\adj(B)\adj(A),
\end{equation}
in combination with \eqref{defGzeta} allows us to simplify \eqref{res1} further:
\begin{equation}
\adj(G_{{S}}(z)-\Id) = \frac{1}{(\ri \tau k)^4}\tilde{\mathbf{Q}}_{\mathbf{k}} \adj(G(\zeta)-\ri \tau k)\tilde{\mathbf{Q}}_{\mathbf{k}}^T.
\end{equation}
For a simple isolated eigenvalue $\lambda$, the kernel of $G$ is one dimensional and hence, there exists a complex function $\zeta\mapsto g(\zeta)$ and a complex vector function $\zeta\mapsto\mathbf{a}(\zeta)$ such that
\begin{equation}\label{tensor}
\adj(G(\zeta)-\ri \tau k ) = g(\zeta)\mathbf{a}(\zeta)\otimes \mathbf{a}^T(\zeta), 
\end{equation}
and consequently,
\begin{equation}
\begin{split}
\adj(G_{{S}}(z)-\Id) & = \frac{1}{(\ri \tau k)^4}\tilde{\mathbf{Q}}_{\mathbf{k}} \adj(G(\zeta)-\ri \tau k)\tilde{\mathbf{Q}}_{\mathbf{k}}^T\\
& = \frac{g(\zeta)}{(\ri \tau k)^4}\tilde{\mathbf{Q}}_{\mathbf{k}} (\mathbf{a}(\zeta)\otimes \mathbf{a}^T(\zeta))\tilde{\mathbf{Q}}_{\mathbf{k}}^T\\
&=\frac{g(\zeta)}{(\ri \tau k)^4} [\tilde{\mathbf{Q}}_{\mathbf{k}}\mathbf{a}(\zeta)\otimes (\tilde{\mathbf{Q}}_{\mathbf{k}}\mathbf{a}(\zeta))^T].
\end{split}
\end{equation}
From \eqref{tensor} it suffices to know one row or column of $\adj(G(\zeta)-\ri \tau k)$ to deduce $\adj(G_{{S}}(\zeta)-\Id)$ completely. Indeed, the last column of $\adj(G(\zeta)-\ri \tau k)$ can be calculated easily and we set 
\begin{equation}\label{a1}
\mathbf{a}(\zeta) = \left(
\begin{array}{c}
\frac{\ri \tau k}{\sqrt{6}}(\zeta+(\zeta^2-1)Z(\zeta))\\
\frac{1}{\sqrt{6}}(1+\ri k\tau\zeta)(\zeta+(\zeta^2-1)Z(\zeta))\\
0\\
0\\
-1-k^2\tau^2-\ri k \tau \zeta -(\ri k \tau+\zeta +\ri k \tau \zeta^2 )Z(\zeta)
\end{array}
\right).
\end{equation}

A lengthy but elementary calculation shows that,
\begin{equation}\label{shearbasis}
    \left.\frac{d}{d\zeta}\adj (G(\zeta)-\ri \tau k)\right|_{\zeta=\ri\frac{\tau\lambda_{\rm shear}+1}{k\tau}} = \begin{pmatrix}
        0 & 0 & 0 & 0 & 0\\
         0 & 0 & 0 & 0 & 0\\
          0 & 0 & A & 0 & 0\\
           0 & 0 & 0 & A & 0\\
            0 & 0 & 0 & 0 & 0\\
    \end{pmatrix},
\end{equation}
for the non-zero complex number
\begin{equation}
\begin{split}
    A=-\frac{i \lambda_{\rm shear}  \left(k^4 \tau^4 +(\lambda_{\rm shear}\tau)^4 +(\tau\lambda_{\rm shear} )^3+\lambda_{\rm shear} \tau^3k^2\right)}{6 k},
    \end{split}
\end{equation}
which gives, again, the two basis vectors \eqref{basisshear} for the eigenspace associated with the shear mode.\\

In conclusion, we see that the approach via Riesz projections is equivalent to direct calculations performed in the previous section. Indeed, dividing the vector in \eqref{a1} by its first entry, we recover the form \eqref{hkexplicit} of the basis vectors for the simple eigenvalues - the exact form of the fifth entry will be determined in the following section. Similarly, we have shown that evaluating the complex residue around the two-fold degenerate eigenvalue $\lambda_{\rm shear}$ in \eqref{shearbasis} produces the same basis vectors as \eqref{basisshear}. In the following section, we will put these basis vectors together to give an explicit description of the coordinate change from spectral variables to macroscopic variables.

\section{Hydrodynamic Equations from Spectral Closure}\label{HydroEqu}

In this section, we derive the evolution equations for the macroscopic variables \eqref{defh} explicitly, based on the change of coordinates \eqref{defhhydroBGK}. First, we analyze the spectral temperature in more detail. As a next step, we describe the transport coefficients arising in the hydrodynamic equations qualitatively and show explicitly how they relate to the eigenvalues.\\

\subsection{Properties of the Spectral Temperature}\label{spectraltemperature}
\begin{figure}
    \centering
    \includegraphics[width=0.6\textwidth]{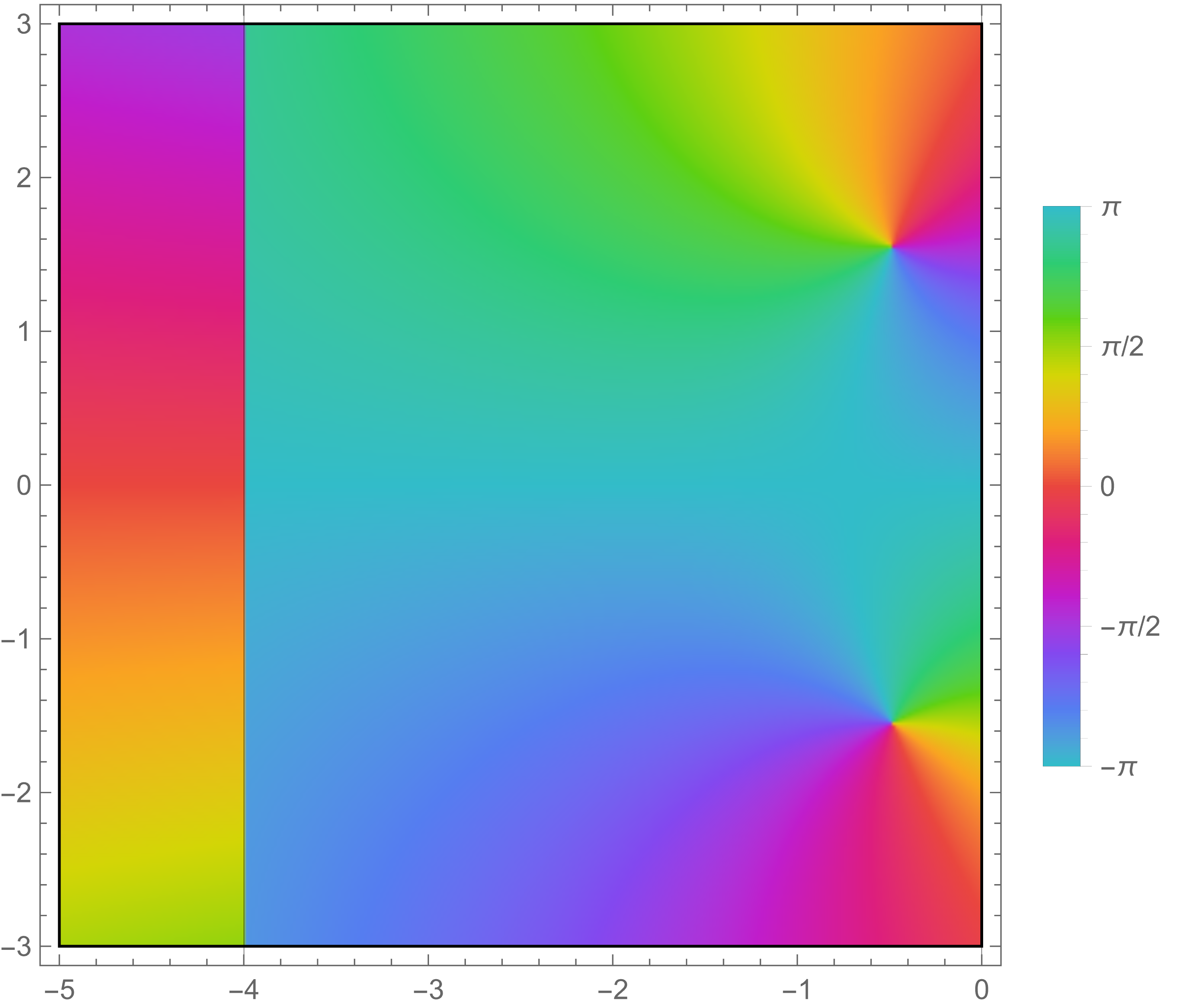}
    \caption{Argument plot of the spectral temperate for $k=1.5$ and $\tau=0.5$, showing two complex conjugate zeros from the numerator in \eqref{spectemp}.}
    \label{fig_Sigma}
\end{figure}

In this subsection, we derive an explicit expression of the spectral temperature and prove some symmetry properties. To this end, we could either evaluate the quotient \eqref{theta1} or just divide \eqref{a1} by its first entry. Indeed, consistency of the two expressions can be checked easily and we proceed by dividing \eqref{a1} by $\frac{\ri \tau k}{\sqrt{6}}(\zeta+(\zeta^2-1)Z(\zeta))$ to recover the $k$-aligned form \eqref{hkexplicit} with explicit spectral temperature
\begin{equation}\label{spectemp}
\begin{split}
    \theta(\lambda) & = \left.\frac{\ri \sqrt{6} \left(k^2 \tau ^2+\ri \zeta  k \tau +Z(\zeta ) \left(\zeta +\ri \left(\zeta ^2+1\right) k \tau \right)+1\right)}{k \tau  \left(\zeta +\left(\zeta ^2-1\right) Z(\zeta )\right)}\right|_{\zeta = \ri\frac{\tau\lambda+1}{\tau k}}\\
    & = \frac{\sqrt{6} \left(\left(k^2 \tau ^2-\tau \lambda  (\tau \lambda +1)\right) Z\left(\frac{\ri (\tau \lambda +1)}{k \tau }\right)-\ri k \tau  \left(k^2 \tau ^2-\tau \lambda \right)\right)}{\left(k^2 \tau
   ^2+(\tau \lambda +1)^2\right) Z\left(\frac{\ri (\tau \lambda +1)}{k \tau }\right)-\ri k \tau  (\tau \lambda +1)}.
\end{split}
\end{equation}
Function \eqref{spectemp} is an analytic function on the strip $\mathcal{R}_{\frac{1}{\tau}}$, see Figure \ref{fig_Sigma}. 

Using \eqref{g0} from Appendix \ref{propZ}, we find that
\begin{equation}
\begin{split}
Z\left(\ri\frac{\tau\lambda+1}{k\tau}\right)^* & =\left[ \left.\ri\sqrt{\frac{\pi}{2}} e^{-\frac{\zeta^2}{2}}\left[\sign(\Im{\zeta})-\erf\left(\frac{-\ri \zeta}{\sqrt{2}}\right)\right] \right|_{\zeta = \ri\frac{\tau\lambda+1}{k\tau} }\right]^*\\
& =  - \left.\ri\sqrt{\frac{\pi}{2}} e^{-\frac{\zeta^2}{2}}\left[-\sign(\Im{\zeta})-\erf\left(\frac{\ri \zeta}{\sqrt{2}}\right)\right] \right|_{\zeta = -\ri\frac{\tau\lambda^*+1}{k\tau} }\\
& = -\left.\ri\sqrt{\frac{\pi}{2}} e^{-\frac{\zeta^2}{2}}\left[-\sign(-\Im{\zeta})-\erf\left(\frac{-\ri \zeta}{\sqrt{2}}\right)\right] \right|_{\zeta = \ri\frac{\tau\lambda^*+1}{k\tau} }\\
& = - Z\left(\ri\frac{\tau\lambda^*+1}{k\tau}\right),
\end{split}
\end{equation}
which implies that
\begin{equation}
    \begin{split}
        \theta(\lambda)^*& = \left[\frac{\sqrt{6} \left(\left(k^2 \tau ^2-\tau \lambda  (\tau \lambda +1)\right) Z\left(\frac{\ri (\tau \lambda +1)}{k \tau }\right)-\ri k \tau  \left(k^2 \tau ^2-\tau \lambda \right)\right)}{\left(k^2 \tau
   ^2+(\tau \lambda +1)^2\right) Z\left(\frac{\ri (\tau \lambda +1)}{k \tau }\right)-\ri k \tau  (\tau \lambda +1)}\right]^*\\
   &= -\frac{\sqrt{6} \left(\left(k^2 \tau ^2-\tau \lambda^*  (\tau \lambda^* +1)\right) Z\left(\frac{\ri (\tau \lambda^* +1)}{k \tau }\right)+\ri k \tau  \left(k^2 \tau ^2-\tau \lambda^* \right)\right)}{-\left(k^2 \tau
   ^2+(\tau \lambda^* +1)^2\right) Z\left(\frac{\ri (\tau \lambda^* +1)}{k \tau }\right)+\ri k \tau  (\tau \lambda^* +1)}\\
       & = \theta(\lambda^*).
    \end{split}
\end{equation}
In particular, we conclude that $\theta|_{\mathbb{R}}\subseteq \mathbb{R}$. This symmetry property will be useful in the following subsection, where we evaluate the spectral temperature on the simple eigenvalues to determine the change of coordinates \eqref{defH} explicitly. 

\subsection{Hydrodynamics in $k$-space}

\begin{figure}
    \centering
    \includegraphics[width=0.5\textwidth]{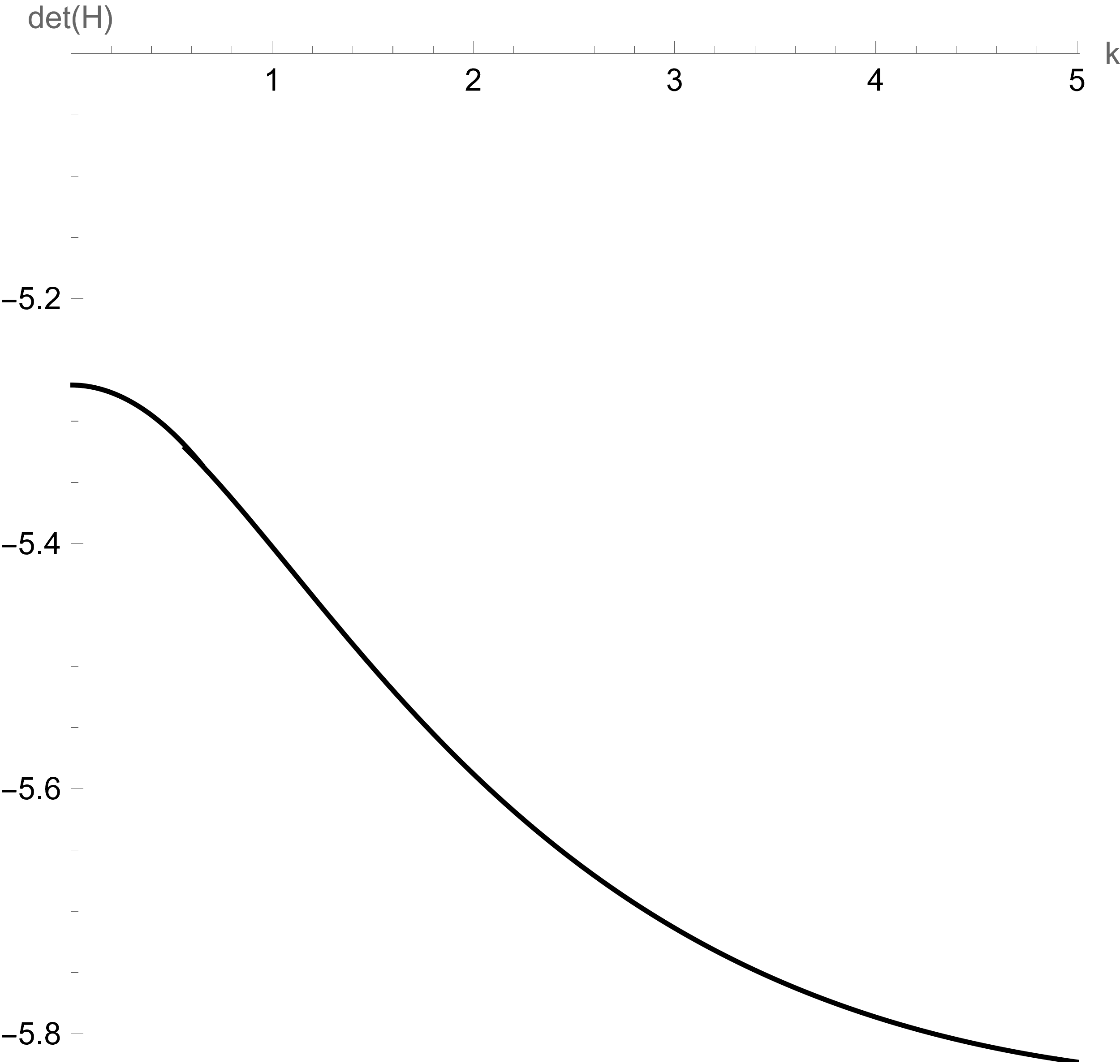}
    \caption{The determinant of the change of coordinates $\mathbf{H}$ for $0\leq k\leq k_{\rm crit,min}$ at $\tau=0.25$.}
    \label{plot_det}
\end{figure}

Using the $k$-aligned basis vectors \eqref{hkexplicit} and \eqref{basisshear}, the change of coordinates from spectral to macroscopic variables \eqref{defH} takes the form,
 \begin{equation}\label{defHfinal}
     \mathbf{H} = \tilde{\mathbf{Q}}_{\mathbf{k}} \begin{pmatrix}
         1 & 1 & 1 & 0 & 0\\
         \frac{\ri }{k}\lambda_{\rm diff} &  \frac{\ri}{k}\lambda_{\rm ac} &  \frac{\ri}{k}\lambda_{\rm ac}^* & 0 & 0\\
         0 & 0 & 0 & 1 & 0\\
         0 & 0 & 0 & 0 & 1\\
         \theta(\lambda_{\rm diff}) &  \theta(\lambda_{\rm ac}) &  \theta(\lambda_{\rm ac}^*) & 0 & 0
     \end{pmatrix}.
 \end{equation}
Its determinant is given by 
 \begin{equation}
 \begin{split}
     \det\mathbf{H} & = \frac{\ri}{k}[( \lambda_{\rm diff}-\lambda_{\rm ac}^*) \theta(\lambda_{\rm ac}) + (\lambda_{\rm ac} - \lambda_{\rm diff})\theta(\lambda_{\rm ac}^*) + (\lambda_{\rm ac}^*-\lambda_{\rm ac}) \theta(\lambda_{\rm diff})]\\
     &=\frac{\ri}{k}\Big(-2\ri(\Im\lambda_{\rm ac})\theta(\lambda_{\rm diff})+2\ri\Im[(\lambda_{\rm ac} - \lambda_{\rm diff})\theta(\lambda_{\rm ac}^*)]\Big)\\
     & = \frac{2}{k}\Big((\Im\lambda_{\rm ac})\theta(\lambda_{\rm diff})-\Im[(\lambda_{\rm ac} - \lambda_{\rm diff})\theta(\lambda_{\rm ac}^*)]\Big),
     \end{split}
 \end{equation}
 which defines a real-valued function of wave number. A plot $k\mapsto \det\mathbf{H}(k)$ is shown in Figure \ref{plot_det}, which already indicates the $\mathbf{H}$ is invertible for all wave numbers $0\leq k\leq k_{\rm crit,min}$. 
Using \eqref{dynalpha} and the invertibility of $\mathbf{H}$, the dynamics for the macroscopic variables on the hydrodynamic manifold are then given by
\begin{equation}\label{hydrodynamics}
\frac{\partial \hat{\mathbf{h}}_{\rm hydro}}{\partial t} = \mathbf{H}\boldsymbol{\Lambda} \mathbf{H}^{-1}\hat{\mathbf{h}}_{\rm hydro}.
\end{equation}

\begin{remark}
We remark that the change of coordinates \eqref{defHfinal} on the hydrodynamic manifold does not involve any expression depending on the shear mode $\lambda_{\rm shear}$ explicitly (see \eqref{basisshear}). 
The \textit{closure operator}, however, will inevitably involve terms that depend on $\lambda_{\rm shear}$ as well. In particular, the hydrodynamics \eqref{hydrodynamics} depend on the shear mode through $\boldsymbol{\Lambda}$. 
\end{remark}

\begin{remark}
We note that the hydrodynamics \eqref{hydro} could be extended beyond the minimal critical wave number by setting $\lambda_{N}(k) =-{1}/{\tau}$ for $k>k_{\rm crit,N}$. Even though, strictly speaking, the eigenvalue does not exist beyond that point, we can, nonetheless, define the hydrodynamic equations by requiring the decay rate of a mode to coincide with the overall minimal decay rate $-{1}/{\tau}$.
\end{remark}

\begin{figure}
     \centering
     \begin{subfigure}[b]{0.3\textwidth}\label{fig_c1}
         \centering
         \includegraphics[width=\textwidth]{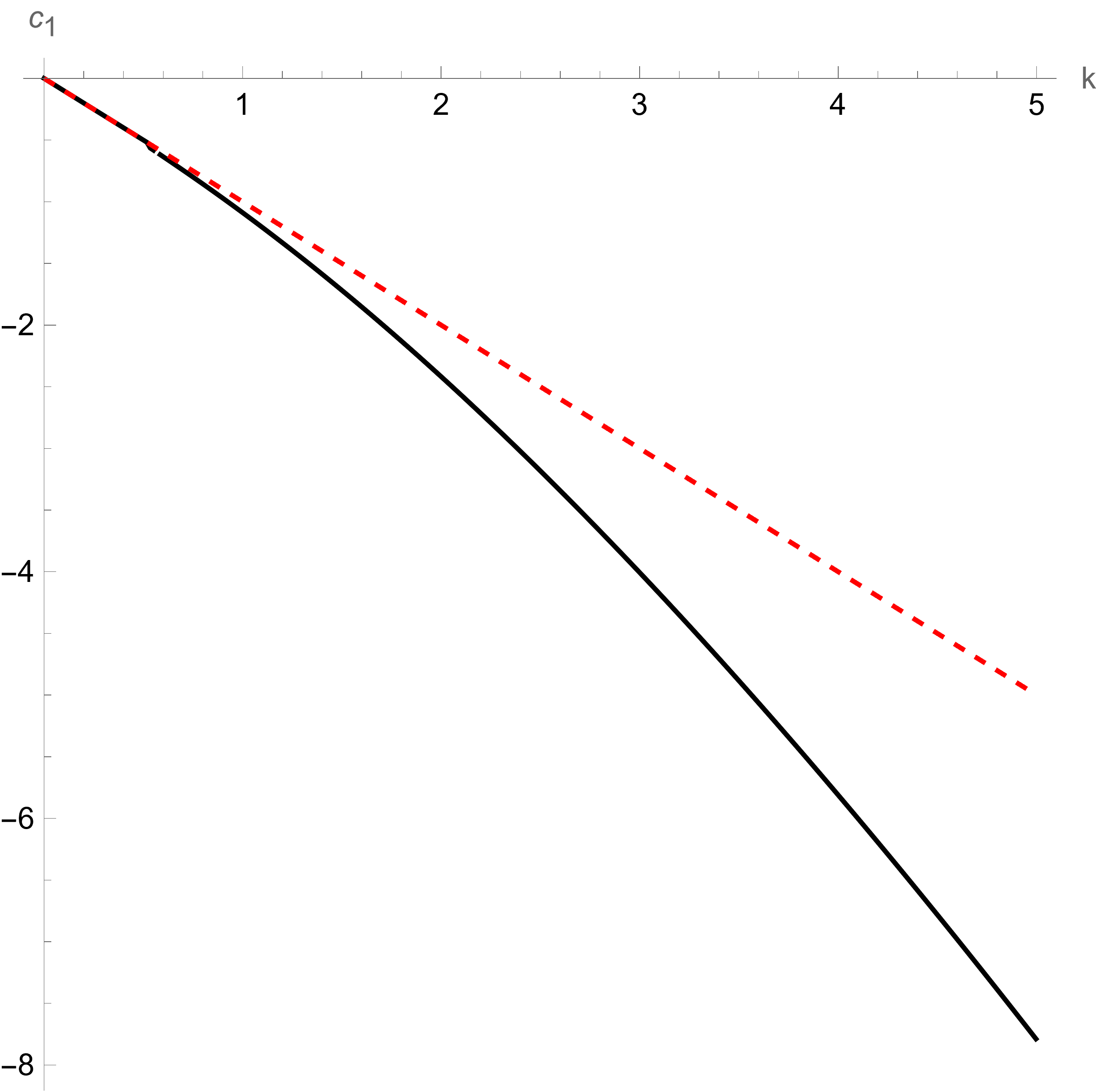}
         \caption{Coefficient $c_1$ in dependence on $k$.}
     \end{subfigure}
     \hfill
     \begin{subfigure}[b]{0.3\textwidth}\label{fig_c2}
         \centering
         \includegraphics[width=\textwidth]{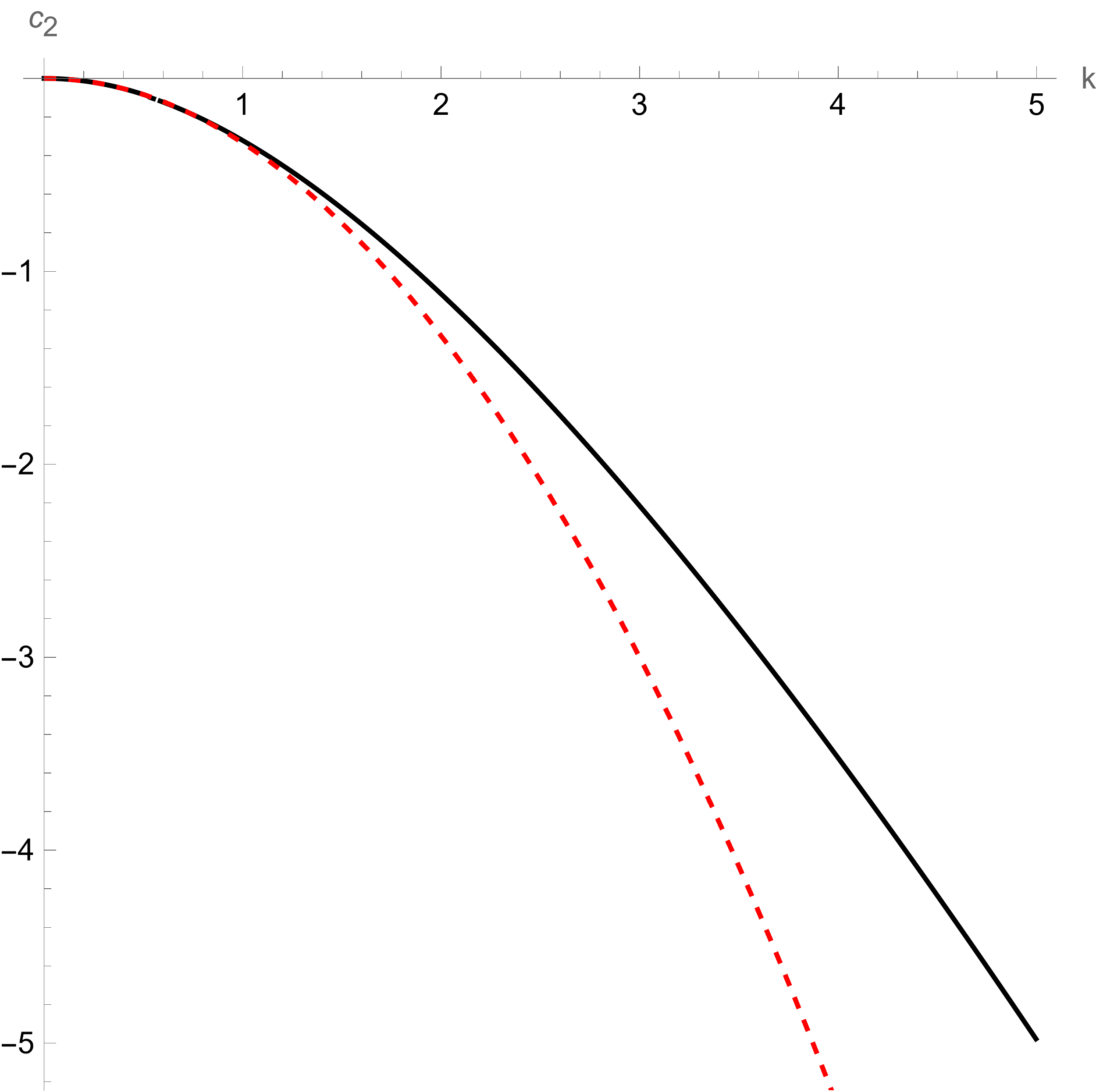}
         \caption{Coefficient $c_2$ in dependence on $k$.}
     \end{subfigure}
     \begin{subfigure}[b]{0.3\textwidth}\label{fig_c3}
         \centering
         \includegraphics[width=\textwidth]{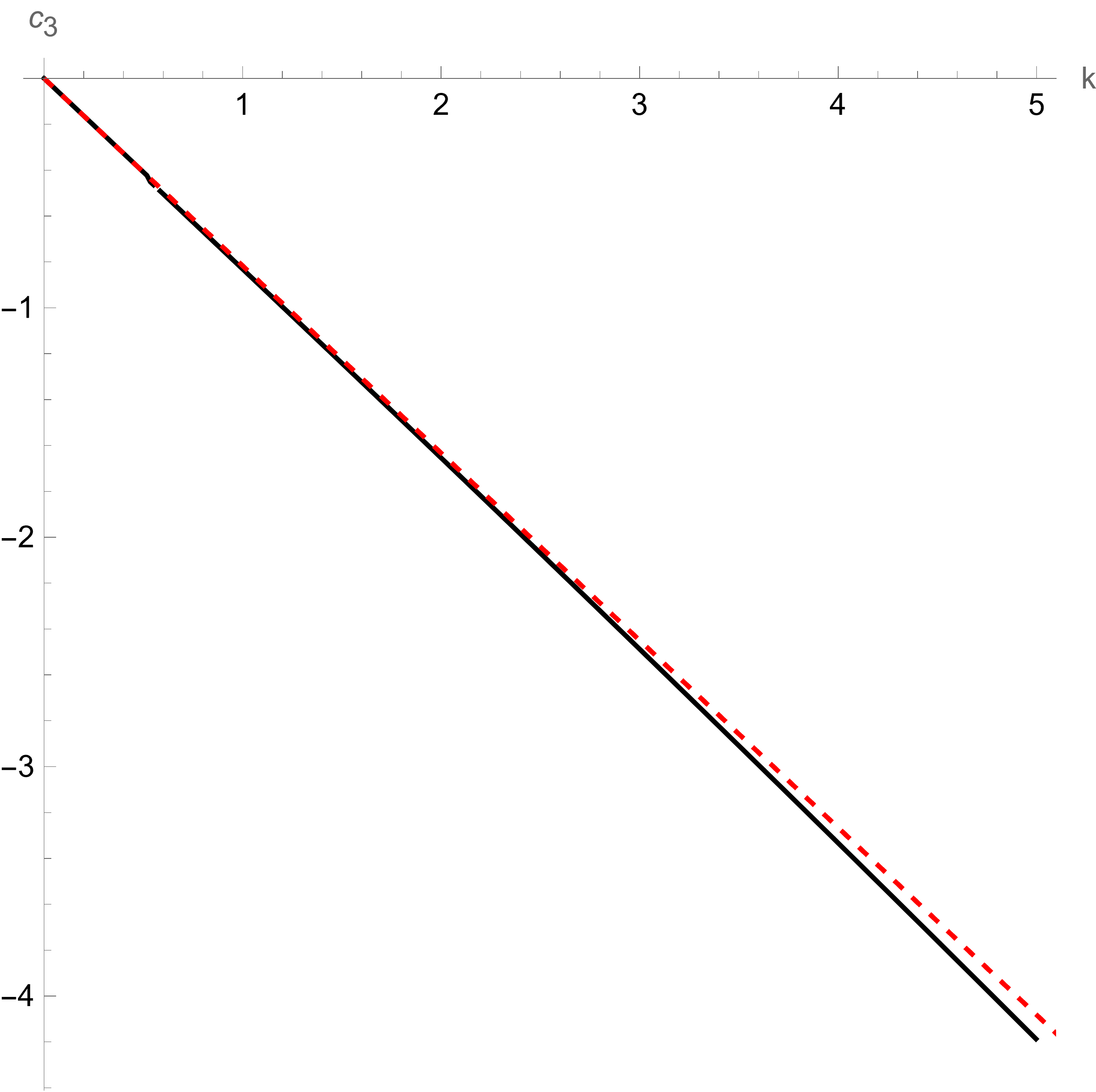}
         \caption{Coefficient $c_3$ in dependence on $k$.}
     \end{subfigure}  
     \begin{subfigure}[b]{0.3\textwidth}\label{fig_c4}
         \centering
         \includegraphics[width=\textwidth]{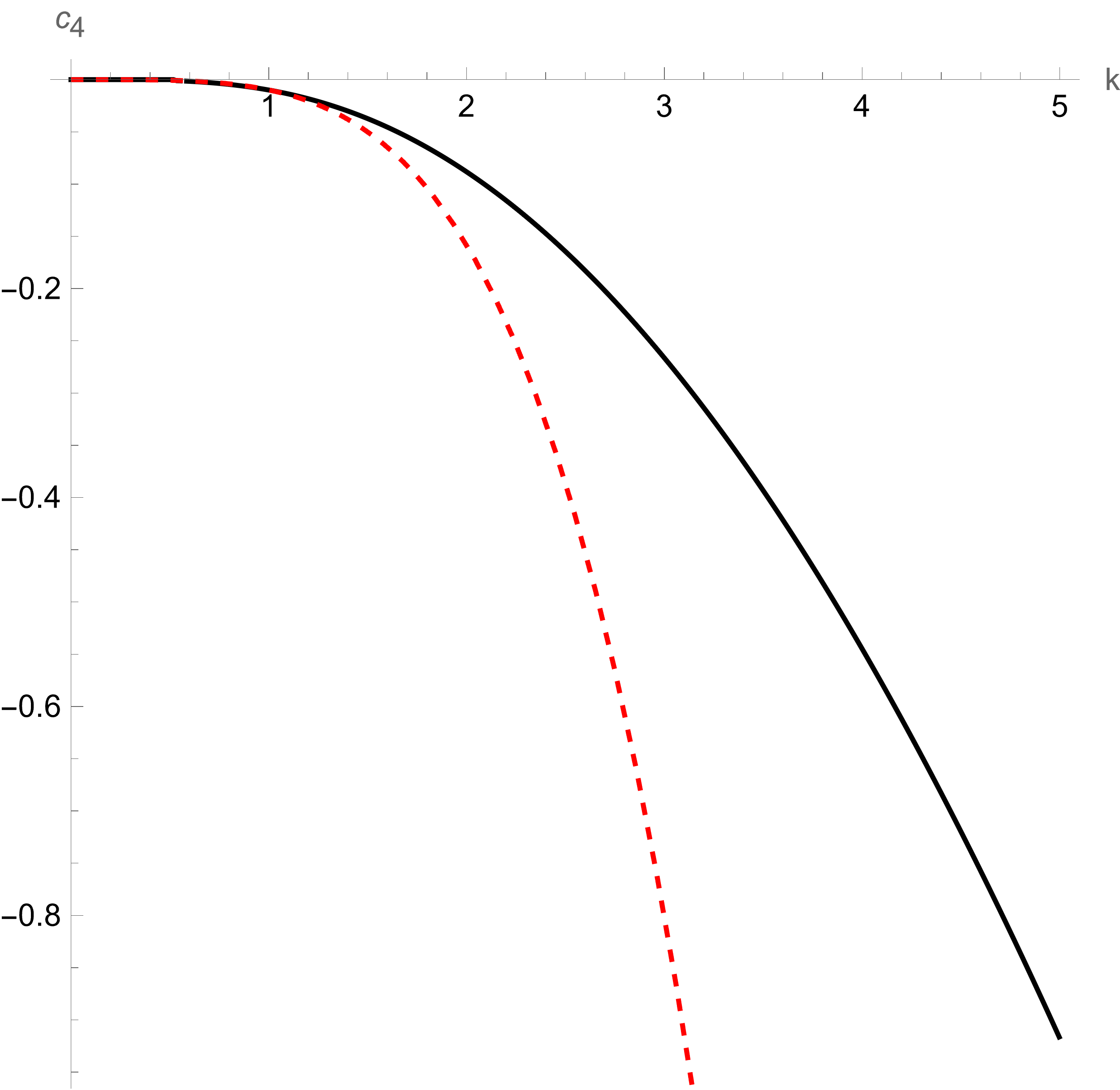}
         \caption{Coefficient $c_4$ in dependence on $k$.}
     \end{subfigure}
     \hfill
     \begin{subfigure}[b]{0.3\textwidth}\label{fig_c5}
         \centering
         \includegraphics[width=\textwidth]{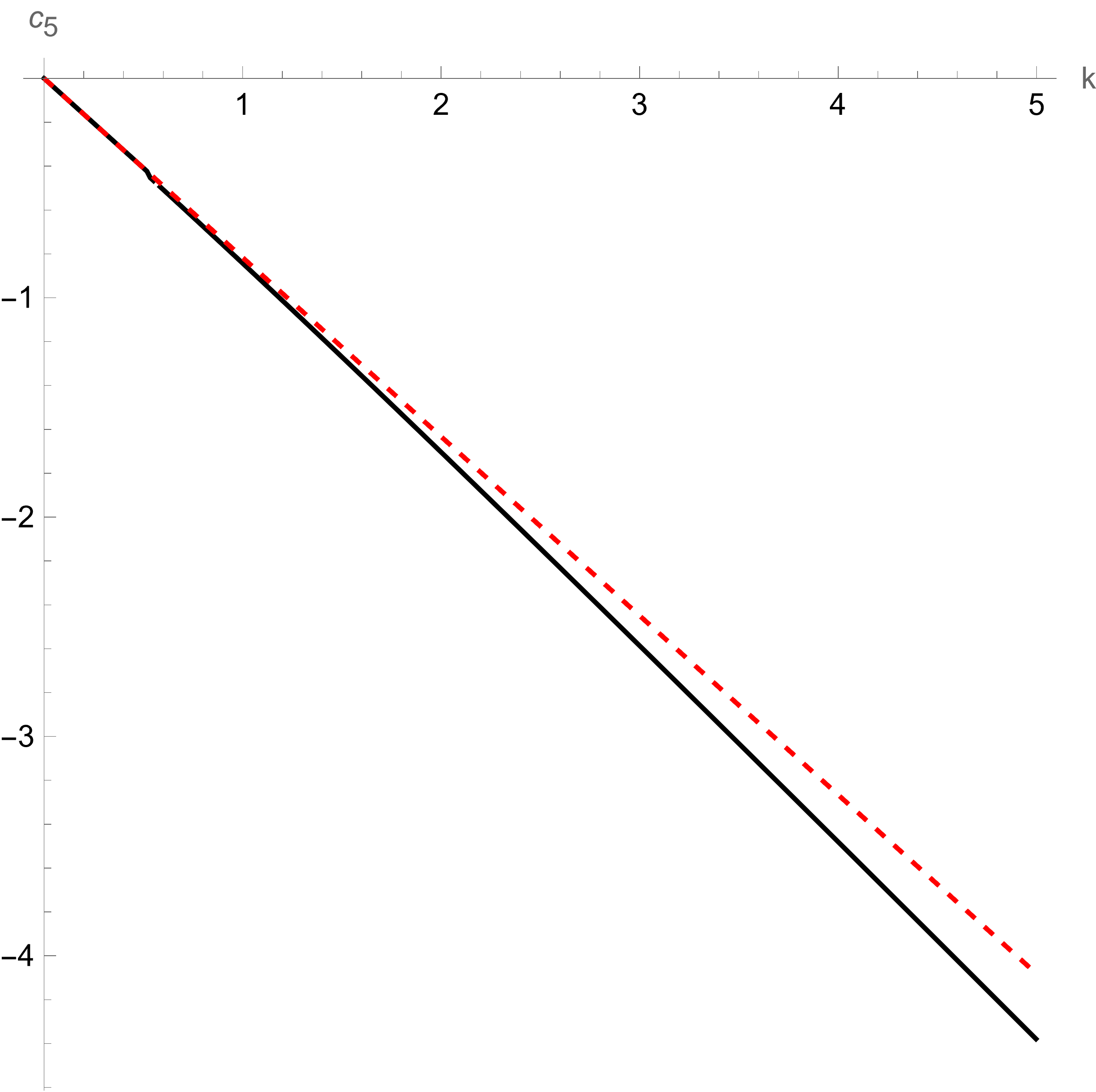}
         \caption{Coefficient $c_5$ in dependence on $k$.}
     \end{subfigure}
     \begin{subfigure}[b]{0.3\textwidth}\label{fig_c6}
         \centering
         \includegraphics[width=\textwidth]{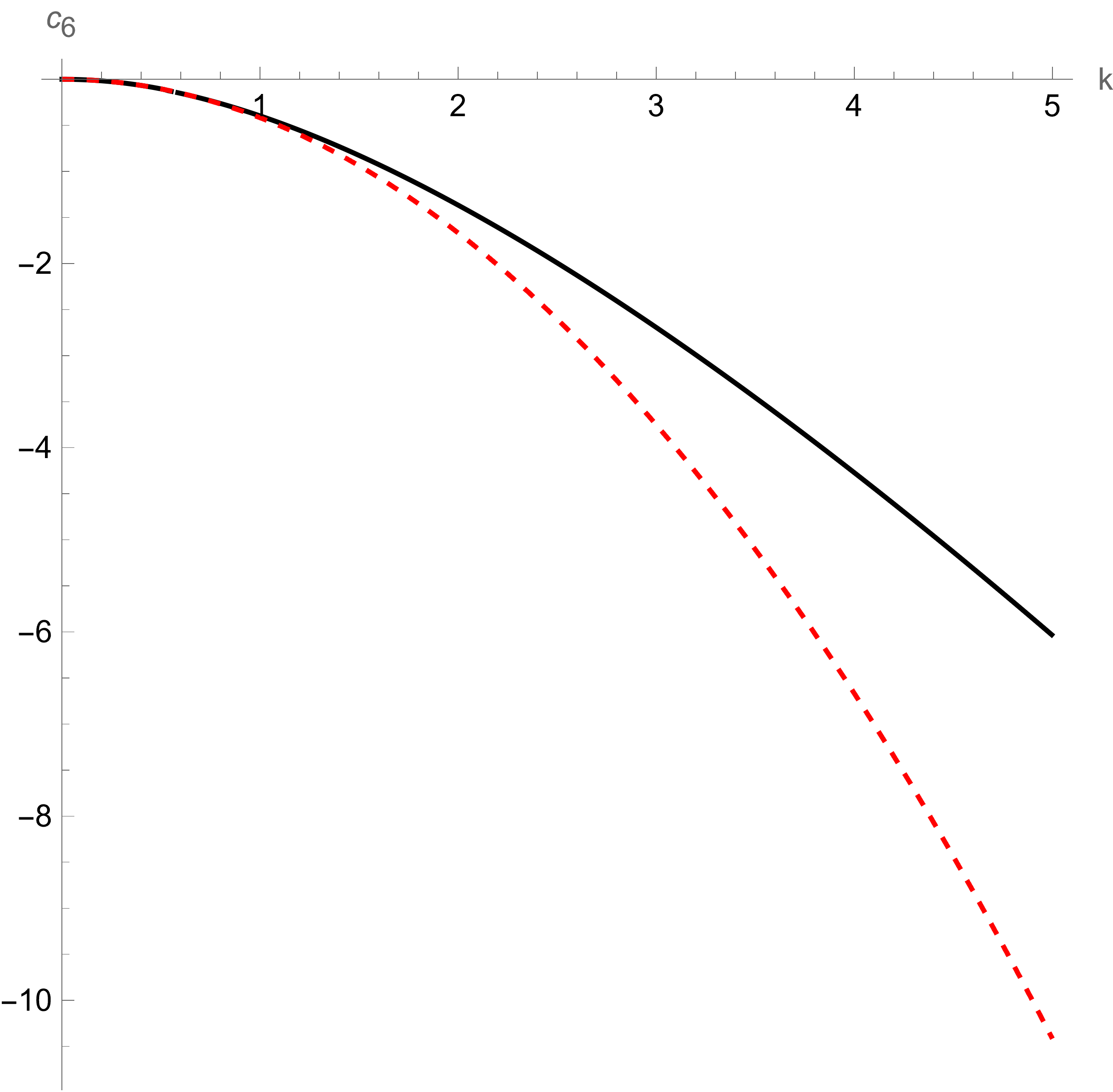}
         \caption{Coefficient $c_6$ in dependence on $k$.}
     \end{subfigure}
     \caption{The transport coefficients calculated in Section \ref{cexpl} in dependence on wave number ($0\leq k\leq k_{\rm crit,min}$) for $\tau=0.25$ (solid black line) compared to its leading-order approximation at the origin (Navier--Stokes/Euler, dashed red line).}
     \label{cplot}
\end{figure}

A cumbersome but elementary calculation shows that
\begin{equation}
    \mathbf{H}\boldsymbol{\Lambda} \mathbf{H}^{-1} = \tilde{\mathbf{Q}}_{\mathbf{k}}  \begin{pmatrix}
         0 & -\ri k & 0 & 0 & 0\\
         C_1 & C_2 & 0 &  0 & C_3\\
         0 & 0 &  \lambda_{\rm shear} & 0 & 0 \\
         0 & 0 & 0 & \lambda_{\rm shear} & 0\\
         C_4 & C_5 & 0 & 0 & C_6
     \end{pmatrix}\tilde{\mathbf{Q}}_{\mathbf{k}}^T,
\end{equation}
for the following cyclic quantities
 \begin{equation}\label{defC}
 \begin{split}
     C_1 & = \frac{1}{k^2\det\mathbf{H}}\sum_{(\lambda_1,\lambda_2,\lambda_3)\in \circlearrowright\boldsymbol{\lambda}_{\rm simple}}\lambda_1\lambda_3(\lambda_1-\lambda_3)\theta(\lambda_2),\\
     C_2 & = \frac{\ri}{k\det\mathbf{H}}\sum_{(\lambda_1,\lambda_2,\lambda_3)\in \circlearrowright\boldsymbol{\lambda}_{\rm simple}} (\lambda_1^2-\lambda_3^2)\theta(\lambda_2),\\
     C_3 &= -\frac{1}{k^2\det\mathbf{H}} \prod_{(\lambda_1,\lambda_2,\lambda_3)\in \circlearrowright\boldsymbol{\lambda}_{\rm simple}} (\lambda_1-\lambda_2),\\
     C_4 & = \frac{\ri}{k\det\mathbf{H}}\sum_{(\lambda_1,\lambda_2,\lambda_3)\in \circlearrowright\boldsymbol{\lambda}_{\rm simple}} \lambda_2 (\lambda_1-\lambda_3)\theta(\lambda_1)\theta(\lambda_3),\\
     C_5 & = \frac{1}{\det\mathbf{H}}\sum_{(\lambda_1,\lambda_2,\lambda_3)\in \circlearrowright\boldsymbol{\lambda}_{\rm simple}}(\lambda_3 - \lambda_1)\theta(\lambda_1)\theta(\lambda_3),\\
     C_6 & = -\frac{\ri}{k\det\mathbf{H}}\sum_{(\lambda_1,\lambda_2,\lambda_3)\in \circlearrowright\boldsymbol{\lambda}_{\rm simple}}\lambda_1\theta(\lambda_1)(\lambda_2-\lambda_3),
     \end{split}
 \end{equation}
 where we have used the notation for cyclical permutations outlined in \eqref{defcirc}.
In Appendix \ref{cexpl}, an explicit expansion of the quantities in \eqref{defC} is performed. We can show that $C_1,C_3$ and $C_4$ are purely imaginary numbers, whereas $C_2,C_5$ and $C_6$  are purely real numbers. Thus we set 
\begin{equation}\label{c}
    C_1 = \ri c_1,\quad C_2 = c_2,\quad C_3 = \ri c_3,\quad C_4 = c_4,\quad C_5 = \ri c_5,\quad C_6 = c_6,
\end{equation}
for $c_j\in\mathbb{R}$, $1\leq j\leq 6$, and the full hydrodynamics become
 \begin{equation}\label{hydro}
     \frac{\partial \hat{\mathbf{h}}_{\rm hydro}}{\partial t} =\tilde{\mathbf{Q}}_{\mathbf{k}}  \begin{pmatrix}
         0 & -\ri k & 0 & 0 & 0\\
         \ri c_1 & c_2 & 0 &  0 & \ri c_3\\
         0 & 0 &  \lambda_{\rm shear} & 0 & 0 \\
         0 & 0 & 0 & \lambda_{\rm shear} & 0\\
          c_4 & \ri c_5 & 0 & 0 & c_6
     \end{pmatrix}\tilde{\mathbf{Q}}_{\mathbf{k}}^T\hat{\mathbf{h}}_{\rm hydro}.
 \end{equation}
Figure \ref{cplot} depicts the coefficients \eqref{c} in dependence of wave number (compared to the Navier--Stokes--Fourier approximation, see Section \ref{compare}).\\

In summary, the hydrodynamic equations in $k$-space are explicitly related to the spectral problem for the linear part through the six transport coefficients $\{c_j\}_{1\leq j\leq 6}$. 

\subsection{Hydrodynamics in real space}

Let us transform equation \eqref{hydro} back to physical coordinates. To this end, we note that
\begin{equation}
\begin{split}
    \mathbf{Q}_{\mathbf{k}}\text{diag}(C_2,\lambda_{\rm shear},\lambda_{\rm shear})\mathbf{Q}_{\mathbf{k}}^T & = \lambda_{\rm shear}\Id_{3\times 3} + \mathbf{Q}_{\mathbf{k}}\text{diag}(C_2-\lambda_{\rm shear},\lambda_{\rm shear},\lambda_{\rm shear})\mathbf{Q}_{\mathbf{k}}^T\\
    &=\lambda_{\rm shear}\Id_{3\times 3} + \frac{1}{k^2}(C_2-\lambda_{\rm shear})\mathbf{k}\otimes\mathbf{k}^T,
    \end{split}
\end{equation}
where we have used the definition of $\mathbf{Q}_{\mathbf{k}}$ in \eqref{defQ}. Consequently, in physical space, the right-hand side of equation \eqref{hydro} translates to a linear integral operator, which can be written as 
\begin{equation}\label{hydrox}
   \frac{\partial }{\partial t}\left(\begin{array}{c}
    \rho\\[0.2cm]
    \mathbf{u}\\[0.2cm]
    T \end{array}\right) = \left(\begin{array}{c}
    -\nabla\cdot\mathbf{u}\\[0.2cm]
     I_1(\Delta)\nabla\rho + I_{\rm shear}(\Delta)\mathbf{u}+I_{2}(\Delta)\nabla(\nabla\cdot\mathbf{u}) + I_3(\Delta)\nabla T \\[0.2cm]
    I_4(\Delta)\rho + I_5(\Delta)\nabla\cdot\mathbf{u} + I_6(\Delta) T
    \end{array}\right),
 \end{equation}
 where the integral operators $\{I_j\}_{1\leq j\leq 6}$ and $I_{\rm shear}$ are related to \eqref{defC} via Fourier series, multiplication/division by $\sqrt{\frac{3}{2}}$ and the rotation matrix \eqref{defQtilde}. The differential operates are to be understood with respect to $\mathbf{x}$ and $\Delta=\nabla\cdot\nabla$ is the Laplacian. We also have used symmetry properties $k\mapsto -k$ from the explicit form in Appendix \ref{cexpl} and the symmetry of eigenvalues.\\

\begin{remark}
Let us summarize the derivation of the non-local, exact hydrodynamics \eqref{hydrox}. The construction begins the evaluation of the eigenvalues by finding zeros of the spectral function \eqref{defSigma}. Since they are given explicitly as the solutions of a transcendental equation, we can analyse them use them in the further analysis. Indeed, up to the minimal critical wave number, there always exists a five-dimensional invariant plane spanned by the eigenvectors of $\mathcal{L}_{\mathbf{k}}$. Using these eigenfunctions, we can construct the spectral closure \eqref{closureBGK} and define the coordinate change \eqref{defH} from spectral variables to macroscopic variables. Once this is achieved, we can write down the exact hydrodynamic equations on the hydrodynamic manifold \eqref{hydro}, which will attract all generic trajectories exponentially fast. The properties of the transport coefficients \eqref{c} can then be analysed in detail. 
\end{remark}

\begin{remark}
Since we posed the governing equations \eqref{maineq} on the three-dimensional torus, only finitely-many wave numbers contribute in the hydrodynamic equation \eqref{hydro}. Consequently, the action of the integral operators  $\{I_j\}_{1\leq j\leq 6}$ and $I_{\rm shear}$ can be written the convolution with an integral kernel of the form
\begin{equation}
    K_{j}(\mathbf{x})= \sum_{k=0}^{k_{\rm crit,min}}\hat{K}_{j}(k) e^{\ri\mathbf{x}\cdot\mathbf{k}},
\end{equation}
where $1\leq j\leq 6$ or $j=\rm shear$ and for coefficients $\hat{K}_{j}(k)\in\mathbb {R}$. Because of the symmetry properties of the eigenvalues, we actually have that  $\hat{K}_{j}(k)=\hat{K}_{j}(k^2)$, corresponding to the dependence of $I_j$ on the Laplacian $\Delta$ only. An effective approximation of the coefficient functions $\hat{K}_j$ will be discussed in a forthcoming paper.
\end{remark}
The quantities $c_2,\lambda_{\rm shear}$ are viscosity terms, while the term $c_6$ can be regarded as non-local (wave-number dependent) version of the thermal diffusivity. We may re-introduce units by scaling according to
     \begin{equation}
        \mathbf{x}\mapsto L^{-1}\mathbf{x},\quad k\mapsto L k,\quad  \mathbf{v}\mapsto v^{-1}_{\rm thermal}\mathbf{v},\quad \tau\mapsto t_{\rm thermal}^{-1}\tau_{\rm relax}, 
     \end{equation}
     for a specific length scale $L$, the \textit{thermal velocity} $v_{\rm thermal}$, the \textit{thermal time} $t_{\rm thermal}$ and the \textit{relaxation time} $\tau_{\rm relax}$.      Given the Boltzmann constant $k_{\rm B}\approx 10^{-23} m^2 kg s^{-2} K^{-1}$, a specific particle mass $m$ and a reference temperature $T_0$, the thermal quantities are defined as
     \begin{equation}
         t_{\rm thermal} = L\sqrt{\frac{m}{k_{\rm B} T_0}},\quad v_{\rm thermal} = \sqrt{\frac{k_{\rm B}}{m} T_0}.
     \end{equation}
     We also define the \textit{mean-free path length} as
     \begin{equation}
         l_{\rm mfp} = \tau L = \tau_{\rm relax} v_{\rm thermal}. 
     \end{equation}
     The macroscopic variables re-scale according to
     \begin{equation}
         \rho \mapsto  \rho_{0}^{-1}\rho,\quad \mathbf{u}\mapsto   v_{\rm thermal}^{-1}\mathbf{u},\quad T\mapsto  T_0^{-1}T,
     \end{equation}
     where $\rho_0$ is a reference density and $T_0$ is a reference temperature. We note that the transport coefficients \eqref{c} can be written as
     \begin{equation}\label{ctauk}
     \begin{split}
         c_1(k,\tau)& = k \tilde{c}_1[(\tau k)^2],\quad  c_2(k,\tau) = \tau^{-1} \tilde{c}_2[(\tau k)^2],\quad c_3(k,\tau) = k \tilde{c}_3[(\tau k)^2],\\
         c_4(k,\tau) &= \tau^{-1} \tilde{c}_4[(\tau k)^2],\quad c_5(k,\tau) = k \tilde{c}_5[(\tau k)^2],\quad c_6(k,\tau) = \tau^{-1} \tilde{c}_6[(\tau k)^2],
         \end{split}
     \end{equation}
see also the asymptotic expansions \eqref{expc2} in the next section. 
Finally, the hydrodynamic equations \eqref{hydro} can be cast in the form 
     \begin{equation}
     \begin{split}
         \frac{\partial \rho}{\partial t} & = -\rho_0\nabla\cdot\mathbf{u},\\
          \frac{\partial \mathbf{u}}{\partial t} & = \frac{k_{\rm B}T_0}{m\rho_0}\mathcal{I}_1[l_{\rm mfp}^2\Delta]\nabla\rho + \frac{1}{\tau_{\rm relax}}\mathcal{I}_{\rm shear}(l_{\rm mfp}^2\Delta)\mathbf{u}+\frac{l_{\rm mfp}^2}{\tau_{\rm relax}}\mathcal{I}_{2}(l_{\rm mfp}^2\Delta)\nabla(\nabla\cdot\mathbf{u}) +\frac{k_{\rm B}}{m\tau_{\rm relax}}\mathcal{I}_3(l_{\rm mfp}^2 \Delta)\nabla T  ,\\
          \frac{\partial T}{\partial t} & = \frac{T_0}{\rho_0\tau_{\rm relax}}\mathcal{I}_4[l_{\rm mfp}^2\Delta]\rho+ T_0\mathcal{I}_5[l_{\rm mfp}^2\Delta](\nabla\cdot\mathbf{u})+ \frac{1}{\tau_{\rm relax}} \mathcal{I}_6[l_{\rm mfp}^2\Delta] T,
      \end{split}   
     \end{equation}
where the integral operators $\{\mathcal{I}_j\}_{1\leq j\leq 6}$ and $\mathcal{I}_{\rm shear}$ are defined through Fourier series and \eqref{ctauk}. 

\section{Comparison to Existing Fluid Models: Small Wave-Number Limit}\label{compare}

In this section, we compare the exact hydrodynamic system \eqref{hydro} to fluid models derived from the Chapmann--Enskog expansion. Because of the coupling between the wave number and the relaxation time through eigenvalues and the k-aligned spectral basis \eqref{defHfinal}, the terms of the Chapmann--Enskog series correspond to an expansion in wave number \eqref{expandlambda}. We write
\begin{equation}\label{Taylorlambda}
\lambda(k) = \sum_{n=1}^\infty\lambda_{n} k^n,
\end{equation}
for any of the four modal branches, for the Taylor expansion of a mode in terms of wave number. Invoking \eqref{Taylorlambda} into the the spectral temperature \eqref{spectemp} and using the asymptotic expansion \eqref{Ipsymptotic} (in the limit $k\to 0$), we can expand
\begin{equation}\label{thetaexpand}
    \begin{split}
        \theta(\lambda(k)) & \sim \frac{\sqrt{6} \left(\left(k^2 \tau ^2-\tau \lambda  (\tau \lambda +1)\right) Z\left(\frac{\ri (\tau \lambda +1)}{k \tau }\right)-\ri k \tau  \left(k^2 \tau ^2-\tau \lambda \right)\right)}{\left(k^2 \tau
   ^2+(\tau \lambda +1)^2\right) Z\left(\frac{\ri (\tau \lambda +1)}{k \tau }\right)-\ri k \tau(\tau \lambda +1)}\\
   & \sim \frac{\sqrt{6} \left(\left(k^2 \tau ^2-\tau \lambda  (\tau \lambda +1)\right) \left(-\frac{k\tau}{\ri(\tau\lambda+1)}-\frac{(k\tau)^3}{[\ri(\tau\lambda+1)]^3}+\mathcal{O}(k^5)\right)-\ri k \tau  \left(k^2 \tau ^2-\tau \lambda \right)\right)}{\left(k^2 \tau
   ^2+(\tau \lambda +1)^2\right) \left(-\frac{k\tau}{\ri(\tau\lambda+1)}-\frac{(k\tau)^3}{[\ri(\tau\lambda+1)]^3}+\mathcal{O}(k^5)\right)-\ri k \tau(\tau \lambda +1)}\\
   & \sim -\sqrt{\frac{3}{2}} \left(\lambda _{{1}}^2+1\right)-\sqrt{\frac{3}{2}} k \lambda _{{1}} \left(\tau  \lambda
   _{{1}}^2+2 \lambda _{{2}}+3 \tau \right)\\
   &\quad-\sqrt{\frac{3}{2}} k^2 \left(3 \tau  \lambda _{{1}}^2 \left(\lambda _{{2}}+\tau \right)+2 \lambda _{{1}} \lambda _{{3}}+3 \tau  \lambda _{{2}}+\lambda _{{2}}^2+3 \tau ^2\right)\\
   &\quad+\sqrt{\frac{3}{2}} k^3 \left(3 \tau ^3 \lambda _{{1}}^3+3 \tau  \lambda _{{1}} \left(-2 \tau  \lambda _{{2}}-\lambda _{{2}}^2+\tau ^2\right)-3 \tau  \lambda _{{1}}^2 \lambda _{{3}}-\lambda _{{3}} \left(2 \lambda _{{2}}+3 \tau\right)\right)\\
   &\quad+\mathcal{O}(k^4),
    \end{split}
\end{equation}
for $k$ sufficiently small. 
Plugging \eqref{thetaexpand} together with \eqref{expandlambda} into \eqref{defHfinal} leads to the following asymptotic expansions for the closure coefficients \eqref{defC}:
\begin{equation}\label{expc1}
    \begin{split}
        C_1 &\sim \frac{i k \left(357 k^6 \tau ^6+991 k^4 \tau ^4-1620 k^2 \tau ^2-900\right)}{60 \left(7 k^2 \tau ^2+15\right)} +\text{h.o.t.},\\
        C_2 & \sim \frac{k^2 \tau  \left(203 k^4 \tau ^4+520 k^2 \tau ^2-600\right)}{30 \left(7 k^2 \tau ^2+15\right)}+\text{h.o.t.}, \\
        C_3 & \sim  -\frac{i k \left(7 k^2 \tau ^2+30\right)^2}{60 \left(7 k^2 \tau ^2+15\right)}+\text{h.o.t.},\\
        C_4 & \sim \frac{k^4 \tau ^3 \left(4437 k^4 \tau ^4-89 k^2 \tau ^2-3000\right)}{90 \left(7 k^2 \tau ^2+15\right)} +\text{h.o.t.},\\
        C_5 & \sim -\frac{i k \left(2523 k^6 \tau ^6+1670 k^4 \tau ^4+360 k^2 \tau ^2+450\right)}{45 \left(7 k^2 \tau ^2+15\right)}+\text{h.o.t.},\\
        C_6 & \sim -\frac{k^2 \tau  \left(203 k^4 \tau ^4+1150 k^2 \tau ^2+750\right)}{30
   \left(7 k^2 \tau ^2+15\right)} +\text{h.o.t.},
    \end{split}
\end{equation}
for small $k$, here $\text{h.o.t.}$ indicates terms of higher order in $k$, either polynomials or rational functions of $k$. Expanding the quotients in \eqref{expc1} in Taylor series around zero, we obtain
\begin{equation}\label{expc2}
    \begin{split}
      C_1 &\sim -\ri k-\ri\frac{4}{3} \tau^2k^3+\mathcal{O}(k^5),\\
      C_2 & \sim -\frac{4}{3}k^2\tau +\frac{16}{9}\tau^3k^4+\mathcal{O}(k^6),\\
      C_3 & \sim -\sqrt{\frac{2}{3}}\ri k+ \mathcal{O}(k^5),\\
      C_4 & \sim -\sqrt{\frac{2}{3}}\frac{10}{3}\tau^3k^4+\mathcal{O}(k^6),\\
      C_5 & \sim -\ri\sqrt{\frac{2}{3}}k-\ri\frac{1}{3}\sqrt{\frac{2}{3}} \tau^2 k^3+\mathcal{O}(k^5),\\
      C_6 & \sim -\frac{5}{3}\tau k^2-\frac{16}{9}\tau^3k^4+\mathcal{O}(k^6), 
    \end{split}
\end{equation}
for small $k$.\\

At first order in $k$, \eqref{expc2}  shows that we have recovered the Euler equation,
\begin{equation}\label{hydroEuler}
     \frac{\partial }{\partial t}  \left(\begin{array}{c}
    \hat{\rho}\\
    \hat{\mathbf{u}}\\
    \hat{T}\end{array}\right)_{\rm Euler}=\tilde{\mathbf{Q}}_{\mathbf{k}}  \begin{pmatrix}
         0 & -\ri k & 0 & 0 & 0\\
         -\ri k & 0 & 0 &  0 &  -\ri k\\
         0 & 0 &  0 & 0 & 0 \\
         0 & 0 & 0 & 0 & 0\\
         0 & -\ri\frac{2}{3} k & 0 & 0 & 0
     \end{pmatrix}\tilde{\mathbf{Q}}_{\mathbf{k}}^T\left(\begin{array}{c}
   \hat{\rho}\\
    \hat{\mathbf{u}}\\
    \hat{T}\end{array}\right)
 \end{equation}
 while at second order in $k$, we recover the Navier--Stokes equation in wave space:
 \begin{equation}\label{hydroNS}
     \frac{\partial }{\partial t}\left(\begin{array}{c}
    \hat{\rho}\\
    \hat{\mathbf{u}}\\
    \hat{T}\end{array}\right)_{\rm Navier-Stokes} =\tilde{\mathbf{Q}}_{\mathbf{k}}  \begin{pmatrix}
         0 & -\ri k & 0 & 0 & 0\\
         -\ri k & -\frac{4}{3}\tau k^2 & 0 &  0 &  -\ri k\\
         0 & 0 &  -\tau k^2 & 0 & 0 \\
         0 & 0 & 0 & -\tau k^2 & 0\\
         0 & -\ri\frac{2}{3} k & 0 & 0 & -\frac{5}{3}\tau k^2
     \end{pmatrix}\tilde{\mathbf{Q}}_{\mathbf{k}}^T\left(\begin{array}{c}
    \hat{\rho}\\
    \hat{\mathbf{u}}\\
    \hat{T}\end{array}\right)
 \end{equation}
Transforming back according to \eqref{defQ} and summing to Fourier series gives the well-known expressions
\begin{equation}\label{hydroEulerx}
     \frac{\partial }{\partial t}  \left(\begin{array}{c}
    \rho\\[0.2cm]
    \mathbf{u}\\[0.2cm]
    T\end{array}\right)_{\rm Euler} = \left(\begin{array}{c}
   -\nabla\cdot \mathbf{u}\\[0.2cm]
   -\nabla (\rho+T)\\[0.2cm]
   -\frac{2}{3}\nabla\cdot \mathbf{u}
    \end{array}
    \right)
 \end{equation}
 as well as 
\begin{equation}\label{hydroNSx}
     \frac{\partial }{\partial t}  \left(\begin{array}{c}
    \rho\\[0.2cm]
    \mathbf{u}\\[0.2cm]
    T\end{array}\right)_{\rm Navier-Stokes} = \left(\begin{array}{c}
   -\nabla\cdot \mathbf{u}\\[0.2cm]
   -\nabla (\rho+T)+\tau\Delta\mathbf{u}+\frac{\tau}{3} \nabla(\nabla\cdot\mathbf{u})\\[0.2cm]
   -\frac{2}{3}\nabla\cdot \mathbf{u}+\frac{5}{3}\tau \Delta T
    \end{array}
    \right)
 \end{equation}
 
Let us comment on the third-order approximation of \eqref{hydro} in $k$, the Burnett equation:
 \begin{equation}\label{Burnett}
     \frac{\partial }{\partial t}\left(\begin{array}{c}
    \hat{\rho}\\
    \hat{\mathbf{u}}\\
    \hat{T}\end{array}\right)_{\rm Burnett} =\tilde{\mathbf{Q}}_{\mathbf{k}}  \begin{pmatrix}
         0 & -\ri k & 0 & 0 & 0\\
         -\ri k -\ri \frac{4}{3} \tau^2k^3 & -\frac{4}{3}\tau k^2 &  &  0 &  -\ri k\\
         0 & 0 &  -\tau k^2 & 0 & 0 \\
         0 & 0 & 0 & -\tau k^2 & 0\\
         0 & -\ri \frac{2}{3}k -\ri \frac{2}{9}\tau^2k^3 & 0 & 0 & -\frac{5}{3}\tau k^2
     \end{pmatrix}\tilde{\mathbf{Q}}_{\mathbf{k}}^T\left(\begin{array}{c}
    \hat{\rho}\\
    \hat{\mathbf{u}}\\
    \hat{T}\end{array}\right)
 \end{equation}
 or, equivalently, in physical space
 \begin{equation}\label{Burnettx}
    \frac{\partial }{\partial t}  \left(\begin{array}{c}
    \rho\\[0.2cm]
    \mathbf{u}\\[0.2cm]
    T\end{array}\right)_{\rm Burnett} = \left(\begin{array}{c}
   -\nabla\cdot \mathbf{u}\\[0.2cm]
   -\nabla (\rho + T)+\frac{4}{3}\tau^2\nabla\Delta\rho+\tau\Delta\mathbf{u}+\frac{\tau}{3} \nabla(\nabla\cdot\mathbf{u})  \\[0.33cm]
   -\frac{2}{3}\nabla\cdot\mathbf{u} - \frac{2}{9}\tau^2\Delta\nabla\cdot\mathbf{u} +\frac{5}{3}\tau  \Delta T
    \end{array}
    \right)
 \end{equation}

 We compare equation \eqref{Burnettx} with  p. 105 equation (34)  in \cite{zheng2004burnett}, the non-dimensional linearized Burnett equation for a more general ellipsoid-statistical BGK (ES-BGK) kinetic model of Holway \cite{holway1965kinetic},
 \begin{equation}\label{Burnettb}
     \begin{split}
        \frac{\partial \rho}{\partial t}&+\nabla\cdot\mathbf{u}=0,\\
        \frac{\partial \mathbf{u}}{\partial t}&+ \nabla(\rho+T)-\Delta\mathbf{u}+\frac{1}{3} \nabla(\nabla\cdot\mathbf{u})-\frac{4}{3}\nabla\Delta\rho-\frac{3b}{3}\nabla\Delta T = 0,\\
        \frac{3}{2}\frac{\partial T}{\partial t}&+\nabla\cdot\mathbf{u}-\frac{5(1-b)}{2}\Delta T-\frac{(1-b)(1-5b)}{3}\Delta\nabla\cdot\mathbf{u},
     \end{split}
 \end{equation}
where parameter $b$ is related to the Prandtl number $\rm Pr = \frac{1}{1-b}$. Since the ES-BGK model reduces to the BGK model with $\rm Pr = 1$ for the BGK equation, we find that equation \eqref{Burnettb} is exactly the same as system \eqref{Burnettx} for $b=0$ and $\tau=1$.
\begin{remark}
As a special feature of the BGK equation, we find that the coefficient $k\mapsto C_{3}(k)$ \eqref{expc2} happens to be very close to the Euler approximation, deviating only at order $k^5$. This explains that there is no contribution of temperature in the non-classical terms in the Burnett approximation, which also implies that the Burnett approximation is globally stable. Indeed, the cubic terms enter only through $C_1$ and $C_5$, which are purely imaginary, and thus only contribute to the higher-order wave motion, not altering the amplitude dynamics. This, of course, is merely a coincidence for the BGK system, since the Shakhov, Maxwell and Hard-Sphere model are not expected to share that property, leading to the Burnett stability. We conjecture that this is due to the temperature coupling back to the velocity dynamics with opposite sign, producing an unstable term.
\end{remark}

\begin{remark}
Since the equations \eqref{hydro} are linear and derived as the invariant dynamics of a globally well-posed system \eqref{maineq} (which is linear itself), the exact hydrodynamics are obviously hyperbolic as a Cauchy problem. The decay rates of solutions on the hydrodynamic manifold, however, will be \textit{weaker} than the decay rate of a general solution, as the hydrodynamic eigenvalues are negative but larger than the general relaxation time in real part. 
\end{remark}

\section{Conclusions and Further Perspectives}\label{sec:conclusion}
We have given an explicit and complete description of the full, non-local hydrodynamic closure of the BGK equation. Based on an explicit description of the spectrum of the linear BGK operator \cite{kogelbauer2023exact}, we obtain an invariant, slow manifold as the space spanned by the hydrodynamic eigenvectors. On this manifold, we are able to explicitly define a closure operator relating the spectral dynamics to the dynamics of the macroscopic variables (density, velocity and temperature) through a linear change of coordinates. The full non-local dynamics are compared to the Euler, the Navier--Stokes--Fourier and the Burnett equations (which may be obtained through the Chapman--Enskog expansion) and full consistency is demonstrated in the small wave-number regime.\\
The explicit form of the transport coefficients in \eqref{hydro} allows us to derive effective approximations in frequency space through polynomials matching both derivatives of the eigenvalues close to zero and the essential spectrum in combination with cut-off functions in wave number. These effective approximations will be non-local as well (involving the convolution with a Dirichlet-type integral kernel), while considerably simplifying the form of the transport coefficients, thus rendering them an interesting candidate for linear gaseous hydrodynamics across all Knudsen numbers.

\section*{Acknowledgement}

This work was supported by European Research Council (ERC) Advanced Grant
834763-PonD. Computational resources at the Swiss National Super Computing
Center CSCS were provided under the grant s1066.

\section*{Declaration of Interest}
The authors declare that there is no conflict of interests.

\appendix 

\section{Explicit Form of the Closure Coefficients}\label{cexpl}
In this section we expand the cyclical expressions for the closure coefficients \eqref{defC} explicitly. This allows us to infer the purely imaginary/real nature in \eqref{c}. 

We expand the first relation in \eqref{defC},
\begin{equation}
    \begin{split}
     C_1 & = \frac{1}{k^2\det\mathbf{H}}\Big[ \lambda_{\rm diff}\lambda_{\rm ac}^*(\lambda_{\rm diff}-\lambda_{\rm ac}^*)\theta(\lambda_{\rm ac})+\lambda_{\rm ac}\lambda_{\rm diff}(\lambda_{\rm ac}-\lambda_{\rm diff})\theta(\lambda_{\rm ac}^*)+\lambda_{\rm ac}^*\lambda_{\rm ac}(\lambda_{\rm ac}^*-\lambda_{\rm ac})\theta(\lambda_{\rm diff})  \Big]\\
&= \frac{2\ri}{k^2\det\mathbf{H}}\Big[ \lambda_{\rm diff}\Im[\lambda_{\rm ac}^*(\lambda_{\rm diff}-\lambda_{\rm ac})\theta(\lambda_{\rm ac})]-|\lambda_{\rm ac}|^2(\Im\lambda_{\rm ac})\theta(\lambda_{\rm diff})  \Big],
    \end{split}
    \end{equation}
 to find that $C_1$ is purely imaginary.\\ 
We expand the second relation in \eqref{defC},
    \begin{equation}
    \begin{split}
     C_2 & = \frac{\ri}{k\det\mathbf{H}}\Big[(\lambda_{\rm diff}^2-(\lambda_{\rm ac}^*)^2)\theta(\lambda_{\rm ac})+(\lambda_{\rm ac}^2-\lambda_{\rm diff}^2)\theta(\lambda_{\rm ac}^*)+((\lambda_{\rm ac}^*)^2-\lambda_{\rm ac}^2)\theta(\lambda_{\rm diff})\Big]\\
     & =  -\frac{2}{k\det\mathbf{H}}\Big[\Im[(\lambda_{\rm diff}^2-(\lambda_{\rm ac}^*)^2)\theta(\lambda_{\rm ac})]-2(\Re\lambda_{\rm ac})(\Im\lambda_{\rm ac})\theta(\lambda_{\rm diff})\Big],
    \end{split}
    \end{equation}
to find that $C_2$ is purely real.\\
We expand the third relation in \eqref{defC},
    \begin{equation}
    \begin{split}
     C_3 & = -\frac{1}{k^2\det\mathbf{H}}(\lambda_{\rm diff}-\lambda_{\rm ac})(\lambda_{\rm ac}-\lambda_{\rm ac}^*)(\lambda_{\rm ac}^*-\lambda_{\rm diff})\\
    & = \frac{2\ri}{k^2\det\mathbf{H}}|\lambda_{\rm diff}-\lambda_{\rm ac}|^2(\Im\lambda_{\rm ac}),
    \end{split}
    \end{equation}
to find that $C_3$ is purely imaginary.\\
We expand the fourth relation in \eqref{defC},
        \begin{equation}
    \begin{split}
     C_4 & = \frac{\ri}{k\det\mathbf{H}} \Big[\lambda_{\rm ac}(\lambda_{\rm diff}-\lambda_{\rm ac}^*)\theta(\lambda_{\rm ac}^*)\theta(\lambda_{\rm diff})+\lambda_{\rm diff}(\lambda_{\rm ac}^*-\lambda_{\rm ac})\theta(\lambda_{\rm ac}^*)\theta(\lambda_{\rm ac})\\
     &\qquad+\lambda_{\rm ac}^*(\lambda_{\rm ac}-\lambda_{\rm diff})\theta(\lambda_{\rm ac})\theta(\lambda_{\rm diff}) \Big]\\
     & = \frac{2}{k\det\mathbf{H}} \Big[\Im[\lambda_{\rm ac}(\lambda_{\rm diff}-\lambda_{\rm ac})\theta(\lambda_{\rm ac})\theta(\lambda_{\rm diff})]+\lambda_{\rm diff}(\Im\lambda_{\rm ac})|\theta(\lambda_{\rm ac})|^2  \Big],
    \end{split}
    \end{equation}
to find that $C_4$ is purely real.\\
We expand the fifth relation in \eqref{defC},
    \begin{equation}
    \begin{split}
     C_5 & = \frac{1}{\det\mathbf{H}} \Big[ (\lambda_{\rm ac}^*-\lambda_{\rm diff})\theta(\lambda_{\rm diff})\theta(\lambda_{\rm ac}^*)+(\lambda_{\rm ac}-\lambda_{\rm ac}^*)\theta(\lambda_{\rm ac})\theta(\lambda_{\rm ac}^*)+(\lambda_{\rm diff}-\lambda_{\rm ac})\theta(\lambda_{\rm diff})\theta(\lambda_{\rm ac}) \Big]\\
    & =\frac{2\ri}{\det\mathbf{H}} \Big[\theta(\lambda_{\rm diff})\Im[\theta(\lambda_{ac})(\lambda_{\rm diff}-\lambda_{\rm ac})]+(\Im\lambda_{\rm ac})|\theta(\lambda_{ac})|^2 \Big],
    \end{split}
    \end{equation}
to find that $c_5$ is purely imaginary.\\
We expand the sixth relation in \eqref{defC},
    \begin{equation}
    \begin{split}
     C_6 & = -\frac{\ri}{k\det\mathbf{H}} \Big[\lambda_{\rm diff}\theta(\lambda_{\rm diff})(\lambda_{\rm ac}-\lambda_{\rm ac}^*)+\lambda_{\rm ac}\theta(\lambda_{\rm ac})(\lambda_{\rm ac}^*-\lambda_{\rm diff})+\lambda_{\rm ac}^*\theta(\lambda_{\rm ac}^*)(\lambda_{\rm diff}-\lambda_{\rm ac})  \Big]\\
    & = \frac{2}{k\det\mathbf{H}} \Big[\lambda_{\rm diff}\theta(\lambda_{\rm diff})(\Im\lambda_{\rm ac})+\Im[\lambda_{\rm ac}\theta(\lambda_{\rm ac})(\lambda_{\rm ac}^*-\lambda_{\rm diff})] \Big],
    \end{split}
    \end{equation}
to find that $C_6$ is purely real.

\section{Properties of the Plasma Dispersion Function $Z$}\label{propZ}

\begin{centering}
\begin{figure}
	\begin{subfigure}{.5\textwidth}
		\centering
		\includegraphics[width=.8\linewidth]{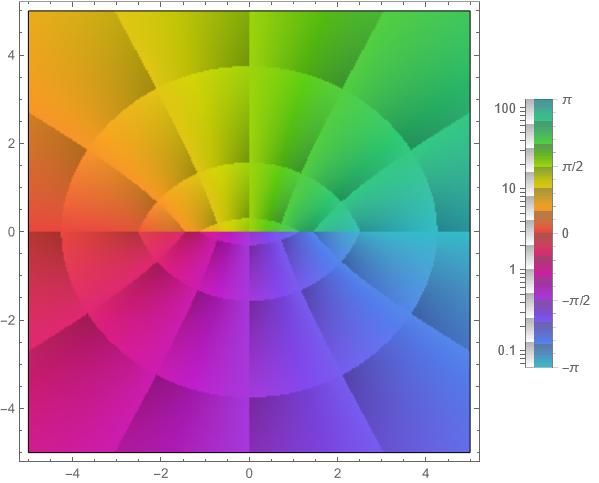}
		\caption{Argument Plot of $Z$}
	\end{subfigure}%
	\begin{subfigure}{.5\textwidth}
		\centering
		\includegraphics[width=1\linewidth]{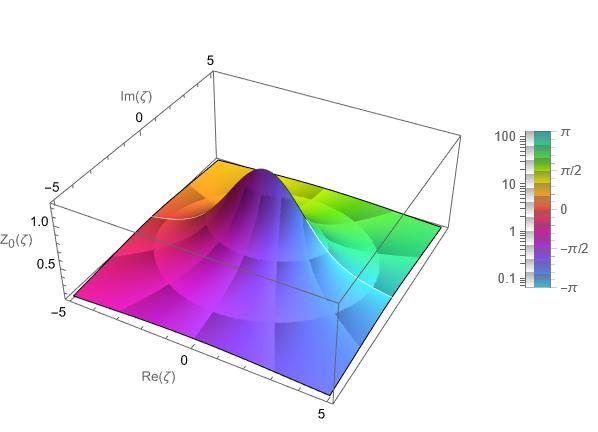}
  \caption{Modulus-Argument Plot of $Z$}
	\end{subfigure}
 \caption{Complex plots of the function Z.}
\label{pI0}
\end{figure}
\end{centering}

In the following, we collect some properties of the plasma dispersion function $Z$, defined through the  integral expression \eqref{defZ}. In our presentation, we will closely follow the calculations performed in \cite{kogelbauer2023exact}.\\
First, let us derive an expression of the integral \eqref{defZ} in terms of less exotic functions. To this end, we rely on the identities in \cite[p.297]{abramowitz1948handbook}. Let
\begin{equation}\label{defw}
w(\zeta)=e^{-\zeta^2}(1-\erf(-\ri \zeta)), \quad \zeta\in\mathbb{C},
\end{equation}
which satisfies the functional identity
\begin{equation}\label{wident}
w(-\zeta)=2e^{-\zeta^2}-w(\zeta),\quad \zeta\in\mathbb{C}.
\end{equation}
Function \eqref{defw} is called \textit{Faddeeva function} and is frequently encountered in problems related to kinetic equations \cite{fitzpatrick2014plasma}. 
We then have that
\begin{equation}
w(\zeta)=\frac{\ri}{\pi}\int_{\mathbb{R}}\frac{e^{-s^2}}{\zeta-s}\, ds,\quad \Im{\zeta}>0,
\end{equation}
and, by relation \eqref{wident}, we have for $\Im{\zeta}<0$:
\begin{equation}
\begin{split}
\frac{\ri}{\pi}\int_{\mathbb{R}}\frac{e^{-s^2}}{\zeta-s}\, ds&=-\frac{\ri}{\pi}\int_{\mathbb{R}}\frac{e^{-s^2}}{(-\zeta)+s}\, ds\\
&=-\frac{\ri}{\pi}\int_{\mathbb{R}}\frac{e^{-s^2}}{(-\zeta)-s}\, ds\\
&=-w(-\zeta)\\
&=e^{-\zeta^2}[-1-\erf(-\ri \zeta)].
\end{split}
\end{equation}
Consequently, we obtain
\begin{equation}
\begin{split}
\int_{\mathbb{R}}\frac{1}{s-\zeta}e^{-\frac{s^2}{2}}\, ds&=\int_{\mathbb{R}}\frac{e^{-s^2}}{s-\frac{\zeta}{\sqrt{2}}}\, ds\\
&=\ri\pi\frac{\ri}{\pi}\int_{\mathbb{R}}\frac{e^{-s^2}}{\frac{\zeta}{\sqrt{2}}-s}\, ds\\
&=\begin{cases}
\ri\pi e^{-\frac{\zeta^2}{2}}\left[1-\erf\left(\frac{-\ri \zeta}{\sqrt{2}}\right)\right],&\quad\text{ if } \Im{\zeta}>0,\\
\ri \pi e^{-\frac{\zeta^2}{2}}\left[-1-\erf\left(\frac{-\ri \zeta}{\sqrt{2}}\right)\right],&\quad\text{ if } \Im{\zeta}<0,
\end{cases}
\end{split}
\end{equation}
where in the first step, we have re-scaled $s\mapsto \sqrt{2}s$ in the integral.
Written more compactly, we arrive at
\begin{equation}\label{g0}
Z(\zeta)=\ri\sqrt{\frac{\pi}{2}} e^{-\frac{\zeta^2}{2}}\left[\sign(\Im{\zeta})-\erf\left(\frac{-\ri \zeta}{\sqrt{2}}\right)\right], \quad \Im{\zeta}\neq 0. 
\end{equation}
An an argument plot together with an modulus-argument plot of $Z$ are shown in Figure \ref{pI0}.
Clearly, $Z$ is discontinuous across the real line (albeit that $Z|_{\mathbb{R}}$ exists in the sense of principal values as the Hilbert transform of a real Gaussian \cite{fitzpatrick2014plasma}). The properties 
\begin{equation}
\begin{split}
    |Z(\zeta)|\leq   \sqrt{\frac{\pi}{2}}, &\text{ for } \zeta\in\mathbb{C}\setminus\mathbb{R},\\
    0< \arg Z(\zeta)<\pi &\text{ for } \Im(\zeta)>0,\\
    -\pi < \arg Z(\zeta) <0 &\text{ for } \Im(\zeta)<0,
\end{split}
\end{equation}
are easy to show and can be read off from the plots \eqref{pI0} directly as well.\\
We also note that
\begin{equation}\label{Z0}
\begin{split}
\lim_{\zeta\to 0,\Im\zeta>0} Z(\zeta) = \ri\sqrt{\frac{\pi}{2}},\\
\lim_{\zeta\to 0,\Im\zeta<0} Z(\zeta) = -\ri\sqrt{\frac{\pi}{2}},
\end{split}
\end{equation}
as can be seen from \eqref{g0}.\\
Function \eqref{g0} satisfies an ordinary differential equation (in the sense of complex analytic functions) on the upper and on the lower half-plane. Indeed, integrating \eqref{defZ} by parts gives
\begin{equation}
\begin{split}
1 &= \frac{1}{\sqrt{2\pi}} \int_{\mathbb{R}}(v-\zeta)\frac{e^{-\frac{v^2}{2}}}{v-\zeta}\, dv=-\zeta Z+\frac{1}{\sqrt{2\pi}}\int_{\mathbb{R}}v\frac{e^{-\frac{v^2}{2}}}{v-\zeta}\, dv\\
&=-\zeta Z-\frac{1}{\sqrt{2\pi}}\int_{\mathbb{R}}\frac{e^{-\frac{v^2}{2}}}{(v-\zeta)^2}\, dv=-\zeta Z-\frac{d}{d\zeta}Z,
\end{split}
\end{equation} 
which implies that $Z$ satisfies the differential equation
\begin{equation}\label{dI0}
\frac{d}{d\zeta}Z= -\zeta Z-1,
\end{equation}
for $\zeta\in\mathbb{C}\setminus\mathbb{R}$. Formula \eqref{dI0} can also be used as a recurrence relation for the higher derivatives of $Z$.\\
Since we will be interested in function \eqref{g0} for $\Im \zeta$ positive and negative as global functions, we define
\begin{equation}\label{defIplus}
\begin{split}
    Z_{+}(\zeta) &= \ri\sqrt{\frac{\pi}{2}} e^{-\frac{\zeta^2}{2}}\left[1-\erf\left(\frac{-\ri \zeta}{\sqrt{2}}\right)\right],\\
    Z_{-}(\zeta) &= \ri\sqrt{\frac{\pi}{2}} e^{-\frac{\zeta^2}{2}}\left[-1-\erf\left(\frac{-\ri \zeta}{\sqrt{2}}\right)\right],
    \end{split}
\end{equation}
for all $\zeta\in\mathbb{C}$. Both functions can be extended to analytic functions on the whole complex plane via analytic continuation.\\
Recall that the error function has the properties that
\begin{equation}
\erf(-\zeta)=-\erf(\zeta),\qquad \erf(\zeta^*)=\erf(\zeta)^*,
\end{equation}
for all $\zeta\in\mathbb{C}$, which implies that for $x\in\mathbb{R}$,
\begin{equation}\label{erfi}
    \erf(\ri x)=-\erf(-\ri x)=-\erf(\ri x)^*,
\end{equation}
i.e, the error function maps imaginary numbers to imaginary numbers. Defining the \textit{imaginary error function},
\begin{equation}
\erfi(\zeta):=-\ri \erf(\ri \zeta),
\end{equation}
for $\zeta\in\mathbb{C}$, which, by \eqref{erfi} satisfies $\erfi|_{\mathbb{R}}\subset\mathbb{R}$, it follows that for $x\in\mathbb{R}$:
\begin{equation}\label{ReImIplus}
  \Re Z_{+}(x)= -\sqrt{\frac{\pi}{2}}e^{-\frac{x^2}{2}}\erfi\left(\frac{x}{\sqrt{2}}\right),\quad \Im Z_{+}(x)= -\sqrt{\frac{\pi}{2}}e^{-\frac{x^2}{2}},
\end{equation}
similarly for $Z_{-}(x)$.\\
Next, let us prove the following asymptotic expansion of $Z_+$:
 \begin{equation}\label{Ipsymptotic}
 Z_{+}(\zeta) \sim -\sum_{n=0}^\infty \frac{(2n-1)!!}{\zeta^{2n+1}},  \qquad \text{ for }|\arg(\zeta)|\leq \frac{\pi}{2}-\delta,\qquad  \zeta\to\infty ,
 \end{equation}
for any  $0<\delta\leq \frac{\pi}{2}$, see also \cite{huba1998nrl}. The proof will be based on a generalized version of Watson's Lemma \cite{https://doi.org/10.1112/plms/s2-17.1.116}. To this end, let us define the Laplace transform
\begin{equation}\label{LaplaceT}
\mathcal{L}[f](\zeta) = \int_0^\infty f(x) e^{-\zeta x}\, dx, \quad \zeta\in\mathbb{C},
\end{equation}
of an integrable function $f: [0,\infty)\to\mathbb{C}$.

\begin{lemma}\label{GWatson}
[Generalized Watson's Lemma]
Assume that \eqref{LaplaceT} exists for some $\zeta=\zeta_0\in\mathbb{C}$ and assume that $f$ admits an asymptotic expansion of the form
\begin{equation}
f(x) =\sum_{n=0}^N a_n x^{\beta_n-1} + o(x^{\beta_N-1}), \qquad x>0,\quad x\to 0,
\end{equation}
where $a_n\in\mathbb{C}$ and $\beta_n\in\mathbb{C}$ with $\Re\beta_0>0$ and $\Re\beta_n>\Re\beta_{n-1}$ for $1\leq n\leq N$.
Then $\mathcal{L}[f](\zeta)$ admits an asymptotic expansion of the form
\begin{equation}
\mathcal{L}[f](\zeta) =\sum_{n=0}^N a_n \Gamma(\beta_n)\zeta^{-\beta_n} +o(\zeta^{-\beta_N}),\quad v,\quad \zeta\to \infty,
\end{equation}
for any real number $0<\delta\leq \frac{\pi}{2}$, where $\Gamma$ is the standard Gamma function. 
\end{lemma}
For a proof of the above Lemma, we refer e.g. to \cite{erdelyi1961general}. Classically, Lemma \eqref{GWatson} is applied to prove that the imaginary error function admits an asymptotic expansion for $x\in\mathbb{R}$ of the form
\begin{equation}\label{erfireal}
    \erfi(x)\sim \frac{e^{x^2}}{\sqrt{\pi}x}\sum_{k=0}^\infty \frac{(2k-1)!!}{(2x^2)^k},\qquad \text{ for } x>0,\quad x\to\infty,
\end{equation}
see also \cite{olver1997asymptotics}, based on the classical version of Watson's Lemma, whose assumptions are, however, unnecessarily restrictive \cite{wong1972generalization}. 

For completeness, we recall the derivation of \eqref{Ipsymptotic} based on Lemma \ref{GWatson}. First, let us rewrite $\erfi$ as a Laplace transform using the change of variables $t=\sqrt{1-s}$ with $dt=\frac{ds}{2\sqrt{1-s}}$
\begin{equation}
    \begin{split}
    \erfi(\zeta)&=\int_0^1\frac{d}{dt} \erfi(t\zeta)\,dt=\frac{2\zeta}{\sqrt{\pi}}\int_0^1  e^{t^2\zeta^2}\, dt = \frac{2\zeta}{\sqrt{\pi}}\int_0^1  e^{\zeta^2(1-s)}\, \frac{ds}{2\sqrt{1-s}}\\
    &= \frac{\zeta e^{\zeta^2}}{\sqrt{\pi}}\int_0^1 \frac{1}{\sqrt{1-s}}  e^{-s\zeta^2}\, ds=\frac{\zeta e^{\zeta^2}}{\sqrt{\pi}}\int_0^{\infty} \frac{\chi_{[0,1]}(s)}{\sqrt{1-s}}  e^{-s\zeta^2}\, ds.
    \end{split}
\end{equation}
From the Taylor expansion of the Binomial function, we know that
\begin{equation}
    \frac{1}{\sqrt{1-s}}=\sum_{n=0}^{\infty}\binom{-1/2}{n} (-s)^n=\sum_{n=0}^{\infty} 4^{-n}\binom{2n}{n}s^n,
\end{equation}
which allows us to apply Lemma \eqref{GWatson} with $\beta_n=n+1$ and $a_n=4^{-n}\binom{2n}{n}$, thus leading to
\begin{equation}
\begin{split}
    \erfi(\zeta)&\sim \frac{\zeta e^{\zeta^2}}{\sqrt{\pi}}\sum_{n=0}^\infty 4^{-n}\binom{2n}{n}\Gamma(n+1) \zeta^{-2(n+1)}\\
    &\sim \frac{e^{\zeta^2}}{\sqrt{\pi}}\sum_{n=0}^\infty \frac{(2n)!}{4^{n}n!} \zeta^{-2n-1}\\
    &\sim \frac{e^{\zeta^2}}{\zeta\sqrt{\pi}}\sum_{n=0}^\infty \frac{(2n-1)!!}{(2\zeta)^{n}},
    \end{split}
\end{equation}
for $\zeta\to\infty$ and $|\arg(\zeta)|\leq \frac{\pi}{2}-\delta$, $0<\delta\leq \frac{\pi}{2}$. This is consistent with formula \eqref{erfireal} for the limit along the real line. Finally, we arrive at the following asymptotic expansion for $Z$:
\begin{equation}
Z_{+}(\zeta)\sim \ri\sqrt{\frac{\pi}{2}}e^{-\frac{\zeta^2}{2}}-\sum_{n=0}^\infty  \frac{(2n-1)!!}{\zeta^{2n+1}},  \qquad \text{ for } |\arg(\zeta)|\leq \frac{\pi}{2}-\delta,\qquad  \zeta\to\infty,
\end{equation}
 which is, of course, equivalent to
 \begin{equation}
 Z_{+}(\zeta) \sim -\sum_{n=0}^\infty \frac{(2n-1)!!}{\zeta^{2n+1}},  \qquad \text{ for }|\arg(\zeta)|\leq \frac{\pi}{2}-\delta,\qquad  \zeta\to\infty ,
 \end{equation}
since $|e^{-\zeta^2}|^2=e^{-2(x^2-y^2)}\to 0$ for $\Re{\zeta}=x\to \infty $.

\bibliographystyle{abbrv}
\bibliography{DynamicalSystems}

\end{document}